%% file: maxsurf-2.tex
\documentclass[12pt]{amsart}
\usepackage{amssymb, epsf, colordvi}
\usepackage{multirow}
\usepackage{latexsym}
\usepackage{graphicx}
\usepackage[all]{xy}

\newif\ifEPSF
\EPSFtrue %
\numberwithin{equation}{section}

\newtheorem{thm}{Theorem}
\numberwithin{thm}{section}

\newtheorem{example}[thm]{Example}
\newtheorem{remark}[thm]{Remark}
\newtheorem{definition}[thm]{Definition}

\newenvironment{rem}{\begin{remark}\rm}{\end{remark}}
\newcounter{FNC}[page]
\def\newfootnote#1{{\addtocounter{FNC}{2}$^\fnsymbol{FNC}$%
     \let\thefootnote\relax\footnotetext{$^\fnsymbol{FNC}$#1}}}

\newcommand\sign{\operatorname{sign}}

\copyrightinfo{}{}

\newcommand{\KK}{\mathbb{K}}

\title{Euler characteristic of primitive $T$-hypersurfaces and maximal surfaces}

\author{Benoit Bertrand}
\address{Section de math\'ematiques\\
Universit\'e de Gen\`eve\\
case postale 64\\
2-4 rue du li\`evre\\
1211 Gen\`eve 4\\
Suisse}
\email{benoit.bertrand@math.unige.ch}
\urladdr{http://name.math.univ-rennes1.fr/benoit.bertrand}

\thanks{The author is grateful to Max Planck Institut f\"ur
  Mathematik in Bonn for excellent working conditions. The author  was partially supported by the European research network IHP-RAAG contract HPRN-CT-2001-00271}

\newtheorem{bbdef}{Definition}[section]

\newtheorem{bbthm}{Theorem}[section]
\newtheorem{bbcor}[bbthm]{Corollary} 
\newtheorem{bblem}{Lemma}

\newtheorem{bbprop}[bbthm]{Proposition}
\newtheorem{bbrmq}{Remark}
\newcommand\CQFD{\hfill $\Box$ \newline}

\newcommand\ZZ{{\mathbb{Z}}}
\newcommand\PP{{\mathbb{P}}}
\newcommand\RR{{\mathbb{R}}}
\newcommand\CC{{\mathbb{C}}}

\newcommand\bin[2]{\binom{#1}{#2}}

\newcommand\bi{\itshape\bfseries}

\newcommand\Star{\operatorname{Star}}
\newcommand\stir{\operatorname{N}}
\newcommand\val{\operatorname{val}}

\newcommand\codim{\operatorname{codim}}
\newcommand\Arg{\operatorname{Arg}}
\newcommand\trop{\operatorname{trop}}

\newcommand\VV{W}
\newcommand\VC{V_\mathbb{C}}

\newcommand{\proofend}{\hfill$\Box$\bigskip}

\begin{document}

\pagestyle{plain}
\setcounter{page}{1}

\begin{abstract}
  Viro method plays an important role in the study of topology of
  real algebraic hypersurfaces. The $T$-primitive hypersurfaces we
  study here appear as the result of Viro's combinatorial patchworking
  when one starts with a primitive triangulation.  We show that the
  Euler characteristic of the real part of such a hypersurface of even
  dimension is equal to the signature of its complex part. We use this
  result to prove the existence of maximal surfaces in some
  three-dimensional toric varieties, namely those corresponding to
  Nakajima polytopes. In fact, these results belong to the field of tropical 
geometry and we explain how they can be understood tropically.
\end{abstract}
\maketitle

\tableofcontents

\section*{Introduction}

The number of connected components of a real algebraic plane
projective curve of degree $m$ is not more than $(m-1)(m-2)/2 + 1$. This
bound was proved by Harnack 
who also showed that for any positive integer $m$ there exist curves
of degree $m$ which are maximal in this sense (i.e.  with
$(m-1)(m-2)/2 + 1$ connected components).  Harnack's bound is
generalized to the case of any real algebraic variety by the
Smith-Thom inequality. Let $b_i(V;K)$ be the $i^{\mbox{\scriptsize
    th}}$ Betti number of a topological space $V$ with coefficients in
a field $K$ (i.e. $b_i(V;K)=\dim_K(H_i(V;K))$).  Denote by $b_*(V;K)$
the sum of the Betti numbers of $V$.  Let $X$ be a complex algebraic
variety equipped with an anti-holomorphic involution $c$. The real
part $\RR X$ of $X$ is the fixed point set of $c$.  Then the
Smith-Thom inequality states that $b_*(\RR X;\ZZ_2) \le b_*(X;\ZZ_2)$.
A variety $X$ for which $b_*(\RR X;\ZZ_2)= b_*(X;\ZZ_2)$ is called a
{\it maximal} variety or {\it $M$-variety}.  The question ``does a
given family of real algebraic varieties contain maximal elements?''
is one of the main problems in topology of real algebraic varieties. For
the family of the hypersurfaces of a given degree in $\RR P^d$ a
positive answer is obtained in \cite{IteVir} using the combinatorial
Viro method called $T$-construction (see \cite{vir3}, \cite{vir6},
\cite{ite1}, and Theorem~\ref{Tcons}).

The Viro patchworking theorem --and its combinatorial particular
case-- is essential in the study of existence of hypersurfaces with
prescribed topology. Let $\Delta$ be a lattice polytope equipped with
a regular (convex) triangulation $\tau$ and a distribution of signs
$D$ at its vertices. From this  data Viro's combinatorial
patchwork produces a real algebraic hypersurface in the toric variety
$X_\Delta$ associated to $\Delta$ whose topology can be easily
recovered from $\Delta, \, \tau $ and $D$ (see the above references
and Section~\ref{tcon}).  {\bi Primitive $T$-hypersurfaces} are those
obtained by the combinatorial patchworking if one starts with a
unimodular (primitive) triangulation $\tau$.  Provided that $\Delta$
correspond to a nonsingular toric variety,
 we prove that for
any primitive $T$-hypersurface $Z$,
\begin{equation}\label{sc}
\chi(\RR Z) = \sigma(Z) \, ,
\end{equation} 
where $\chi(\RR Z)$ is the Euler characteristic of $\RR Z$ and
$\sigma(Z) := \sum_{p+q=0\> [2]} {(-1)}^p h^{p,q}(Z)$ is the signature
of $Z$ if $\dim Z$ is even and $0$ otherwise (See Theorem~\ref{theprop}).

This result can be naturally formulated in the tropical setting introduced
by G. Mikhalkin. (See Section~\ref{S:trophyp} or \cite{Mikhmath.AG/0403015} and
\cite{Mikh05} for definitions and propreties of tropical hypersurfaces.):

The Euler characteristic of a compactified nonsingular real tropical
hypersurface with Newton polytope $\Delta$ is equal to the signature
of a complex algebraic hypersurface in $X_\Delta$ with Newton polytope
$\Delta$.
  
For primitive $T$-surfaces in the projective space $\PP^3$ this result
was first proven by I.~Itenberg \cite{ite1} who used it to construct
maximal $T$-surfaces in all degrees in $\PP^3$. The existence of
maximal surfaces in all degrees in the projective space had been
proven in \cite{Vir79} by Viro using small perturbation technics.

Using the combinatorial Viro method and equality~(\ref{sc}) we prove
the existence of maximal surfaces in a family of toric varieties which
includes $\PP^3$. A polytope $\Delta$ in $\RR ^d$ is a {\bi Nakajima
  polytope} (see Figure \ref{figure:nakapol3}) if either $\Delta$ is
$0$-dimensional or there exists a Nakajima polytope $\bar{\Delta}$ in
$\RR^{d-1}$ and a linear function $f: \RR^{d-1} \to \RR$ , nonnegative
on $\bar{\Delta}$, such that $f(\ZZ ^{d-1})\subset \ZZ$ and
$\Delta=\{(x,x_d) \in \bar{\Delta}\times \RR \mid 0 \le x_d \le
f(x)\}$. We show that for any 3-dimensional Nakajima polytope
$\Delta$ corresponding to a nonsingular toric variety there exists a
maximal hypersurface in $X_\Delta$ with Newton polytope $\Delta$.

I am grateful to Ilia Itenberg for his valuable advice and to Erwan
Brugall\'e for his pertinent remarks.

{\bf Organization of the material.} We first introduce notation and
describe combinatorial patchworking in Section~\ref{prelim}. In the
same section, we describe briefly tropical hypersurfaces and state the
results in the tropical setting. We also recall there facts about
lattice triangulations and state Danilov and Khovanskii Formulae.
Section~\ref{sigmachi} is devoted to the proof of equality~(\ref{sc})
for primitive hypersurfaces stated in Theorem~\ref{theprop}. In the
third section we construct maximal surfaces in a family of toric
varieties using equality~(\ref{sc}). We postpone some proofs of Lemmae
needed in Section~\ref{sigmachi} to the appendix
(Section~\ref{appendix}) where we also give some usefull combinatorial
formulae (e.g. Lemma~\ref{vanlint}).

\section{Preliminaries}\label{prelim}

\subsection{Toric varieties}

We fix here some conventions and notations, the construction of toric
varieties we use is based on the one described in \cite{ful}.  Fix an
orthonormal basis of $\RR ^d$ and thus an inclusion $\ZZ ^d \to
\RR^d$. This inclusion defines a lattice $N$ in $\RR^d$.  Denote by
$M$ the dual lattice $\text{\rm Hom}(N, \ZZ)$ of $N$.  Here we
consider polytopes in $M\otimes\RR$.  By {\bi polytope} we mean
convex polytope whose vertices are integer (i.e. belong to the lattice
$M$).  Let $\Delta$ be a polytope in $M\otimes\RR$.  Let $p$ be a
vertex of $\Delta$ and let $\Gamma_1, \cdots, \Gamma_k$ be the facets
of $\Delta$ containing $p$. To $p$ we associate the cone $\sigma_p$
generated by the 
minimal integer inner normal vectors of $\Gamma_1, \cdots,
\Gamma_k$.  The inner normal fan $\mathfrak{E}_\Delta$ is the fan
whose $d$-dimensional cones are the cones $\sigma_p$ for all vertices
$p$ of $\Delta$. The toric variety $X_\Delta$ associated to $\Delta$ is
the toric variety $X(\mathfrak{E}_\Delta)$ associated to the fan
$\mathfrak{E}_\Delta$.

\subsection{Combinatorial
patchworking}\label{tcon}

Let us  briefly describe the {\it combinatorial patchworking}, also called
{\it  $T$-construc\-tion} , which  is a particular
case  of  the Viro  method. A more detailed exposition can be found in \cite{IteVir}, \cite{vir6} or \cite{GKZ}~p.~385.

By a {\bi subdivision} of a polytope we mean a subdivision in
convex polytopes whose vertices have integer coordinates.  A
subdivision $\tau$ of a $d$-dimensional polytope $\Delta$ is called
{\bi convex} (or {\bi regular}) if there exists a convex
piecewise-linear function $\Phi : \Delta \to \RR$ whose domains of
linearity coincide with the $d$-dimensional polytopes of $\tau$.
A triangulation is a subdivision into simplices.

Given a triple $(\Delta,\tau,D)$, where $\Delta$ is a polytope, $\tau$
a convex triangulation of $\Delta$, and $D$ a distribution of signs at
the vertices of $\tau$ (i.e. each vertex is labelled
with $+$ or $-$), the combinatorial patchworking, produces an
algebraic hypersurface $Z$ in $X_\Delta$.

Let $\Delta$ be a $d$-dimensional polytope in the positive orthant
$(\RR^+)^d = \{(x_1, \ldots x_d) \in \RR^d \mid\\ x_1 \geq 0, \ldots
, x_d \geq 0\}$, and $\tau$ be a convex triangulation of $\Delta$.
Denote by $ s_{(i)}$ the reflection with respect to the coordinate
hyperplane $x_i=0 $ in $\RR^d$.  Consider the union $\Delta^*$ of all
copies of $\Delta$ under the compositions of reflections $ s_{(i)} $
and extend $\tau$ to a triangulation $\tau^*$ of $\Delta^*$ by means
of these reflections.  We extend $D(\tau)$ to a distribution of signs at
the vertices of $\tau^*$ using the following rule : for a vertex $a$
of $\tau^*$, one has $\sign(s_{(i)}(a))=\sign(a)$ if the $i$-th
coordinate of $a$ is even, and $\sign(s_{(i)}(a))=-\sign(a)$,
otherwise.

Let $\sigma$ be
a $d$-dimensional simplex of  $\tau^*$ with vertices of different signs, and
$E$ be the piece of hyperplane which is the convex hull of the middle points of the
edges of $\sigma$
with endpoints of opposite signs. We separate vertices of
$\sigma$ labelled with $+$ from vertices labelled with $-$ by $E$.
The union of all these hyperplane pieces forms a piecewise-linear
hypersurface $ H^*$.

For any facet $ \Gamma$ of $ \Delta^*$, let $N^\Gamma$ be a vector
normal to $\Gamma$. Let $F$ be a face of $\Delta^*$ and $\Gamma_1,
\dots, \Gamma_k$ be the facets containing $F$.  Let $L$ be the linear
space spanned by $N^{\Gamma_1}, \dots, N^{\Gamma_k}$. For any $v =
(v_1, \dots, v_d) \in L \cap \ZZ^d$ identify $F$ with ${s_{(1)}}^{v_1}
\circ {s_{(2)}}^{v_2} \circ \dots \circ {s_{(d)}}^{v_d}(F)$.  Denote
by $ \widetilde{\Delta}$ the result of the identifications.  The
variety $\widetilde{\Delta}$ is homeomorphic to the real part $\RR
X_\Delta$ of $X_\Delta$ (see, for example, \cite{GKZ} Theorem~5.4 p.
383 or \cite{stu} Proposition 2).

Denote by $\widetilde{H}$ the image of $H^*$ in $ \widetilde{\Delta}$.
Let $Q$ be a polynomial with Newton polytope $\Delta$. It defines a
hypersurface $Z_0$ in the torus $(\CC^*)^d$ contained in $X_\Delta$.
The closure $Z$ of $Z_0$ in $X_\Delta$ is the hypersurface defined by
$Q$ in $X_\Delta$. We call $\Delta$ the { \bi Newton polytope } of
$Z$.

\begin{bbthm}[T-construction, O. Viro]\label{Tcons}
Under the  hypotheses made above, there exists a
hypersurface $Z$ of $X_\Delta$ with
Newton polytope $\Delta$
and a homeomorphism $h: \RR X_\Delta\to
\widetilde{\Delta}$ such that  $h(\RR  Z) = \widetilde{H}$.
\end{bbthm}

A  {\bi real algebraic $T$-hypersurface} is a hypersurface which can
be obtained via the above construction.
  A $d$-dimensional simplex with integer
vertices is called {\bi primitive} (or {\bi unimodular}) if its volume
is equal to $\frac{1}{d !}$.  A triangulation $\tau$ of a
$d$-dimensional polytope is {\itshape {\bfseries primitive}} if every
$d$-simplex of the triangulation is primitive. A $T$-hypersurface $Z$
is {\bi primitive} if 
 it can be constructed with the combinatorial patchworking using a
primitive triangulation.

Let $\Delta$ be a $d$-dimensional polytope and $\tau$ a (convex)
triangulation of $\Delta$.  The star, $\Star(v)$, of a vertex $v$ of
$\tau$ is the union of all simplices of $\tau$ containing $v$. Suppose that
$\tau$ is equipped with a distribution of signs. We say that $v$ is
{\bi isolated} if all other vertices of $\Star(v)$ have the sign opposite
to $\sign(v)$.

\begin{bbrmq}\label{sphere}
Let  $p$  be an  isolated  vertex  of  the triangulation  $\tau^*$  of
$\Delta^*$. Assume that $p$ is in the interior of $\Delta^*$. Then the
star of $p$ contains a  connected component of $H^*$ homeomorphic to a
sphere $S^{d-1}$.  So, $p$ corresponds to a $(d-1)$-dimensional sphere of the
real part of the hypersurface $Z$ obtained by Viro's theorem.
\end{bbrmq}

Notice that the triangulation $\tau^*$
induces a natural cell decomposition of $H^*$: each cell is the
intersection of a simplex of $\tau^*$ with $H^*$.
A $k$-dimensional simplex in $\RR^d$ is called {\bi elementary} if the
reductions modulo $2$ of its vertices generate a $k$-dimensional
affine space over $\ZZ_2$.

  Our computation of the Euler characteristic of $\RR Z$ uses the
  following proposition due to Itenberg (see \cite{ite1} Proposition
  3.1).

\begin{bbprop}\label{itenprop}
  Let $s$ be an elementary $k$-simplex of a triangulation $\tau$ of a
  $d$-dimensional polytope.  Assume that $s$ is contained in $j$
  coordinate hyperplanes.  Then the union $s^*$ of the symmetric copies  of
  $s$ contains exactly $(2^d-2^{d-k})/2^j$ cells of dimension $k-1$ of
  the cell decomposition of $H^*$.
\end{bbprop}

\subsection{Tropical hypersurfaces}\label{S:trophyp}

Before stating the results we briefly recall some basic facts about
tropical hypersurfaces. A detailed presentation can be found in
\cite{Mikh05}. In this article we will use the definition of tropical
varieties as images of algebraic varieties defined over the field of
Puiseux series.

Let $\KK$ be the field of Puiseux series. 
 An element of $\KK$ is a
series $g(t)=\sum_{r\in R} b_r t^r$ 
where $b_r \in \mathbb{C}$, $ R \subset
\mathbb{Q}$ is bounded from below and contained in an arithmetic sequence.
Consider the valuation 
 $\val (g(t)):=\min \{r | b_r \neq 0 \}$. For convenience we choose 
Mikhalkin's conventions and consider in fact minus the valuation, 
  $v(g):=-\val(g)$. Let $f$ be a polynomial in $\KK[z_1,\cdots,z_n]
  =\KK[z]$. It is of the form $f(z)=\sum_{\omega \in A} \;c_\omega
  z^\omega$  with $A$ a finite subset of $\ZZ^n$ and $c_\omega \in
  \KK^*$. Consider $Z_f:=\{z\in (\KK ^*)^n | f(z)=0 \}$ the zero set
   of $f$ in  $(\KK ^*)^n$
and the coordinatewise valuation (up to sign)
\[
\begin{array}{rrrl} 
 V: &(\KK ^*)^n &\longrightarrow &\RR^n \\
    & z  &  \longmapsto &(v(z_1), \dots, v(z_n)).
\end{array}
\]

\begin{definition}
  The tropical hypersurface $Z_f^{\trop}$ associated to $f$ is the
  closure (in the usual topology) of the image under $V$ of the
  hypersurface $Z_f$:

$$ Z^{\trop}_f := \overline{V(Z_f)} \subset \RR^n.$$
\end{definition}

The following theorem, due to Kapranov, allows one to see tropical
hypersurfaces as nonlinearity domains of some piecewise-linear convex
functions. Put

\[
\begin{array}{rrrl} 
 \nu: & A &\longrightarrow &\RR \\
                & \omega  &  \longmapsto & - v (c_\omega)\, ,
\end{array}
\] 

the Legendre transform $\mathcal{L}(\nu)$ of $\nu$  is the
piecewise-linear convex function defined as follows:
\[
\begin{array}{rrrl} 
\mathcal{L}(\nu): &\RR^n &\longrightarrow &\RR \\
                & x  &  \longmapsto & \max_{\omega\in A} (x \cdot \omega - \nu(\omega)) 
\end{array}
\]

\begin{thm}[Kapranov]
The tropical hypersurface $Z^{\trop}_f$ is the nonlinearity domain of  $\mathcal{L}(\nu)$:

$$\overline{V(Z_f)}=\mbox{corner locus}(x\mapsto \max_{\omega\in A} (x \cdot \omega +
  v(c_\omega)) ).$$ 
\end{thm} 
 
Another very nice fact is that one can naturally associate to a
tropical hypersurface a dual subdivision of its Newton polytope.
Let $ \Gamma$ be the convex hull in $\RR^n\times\RR$ of the points
$(\omega,v(c_\omega))$ for all $\omega$ in $A$. Put
$$
\begin{array}{rrrl} 
\overline{\nu} : & \Delta &\longrightarrow &\RR \\
                & x &  \longmapsto & \min \{y | (x,y) \in \Gamma\}.
\end{array}
$$

\begin{figure}
\begin{tabular}{ccc}
\includegraphics[width=5cm]{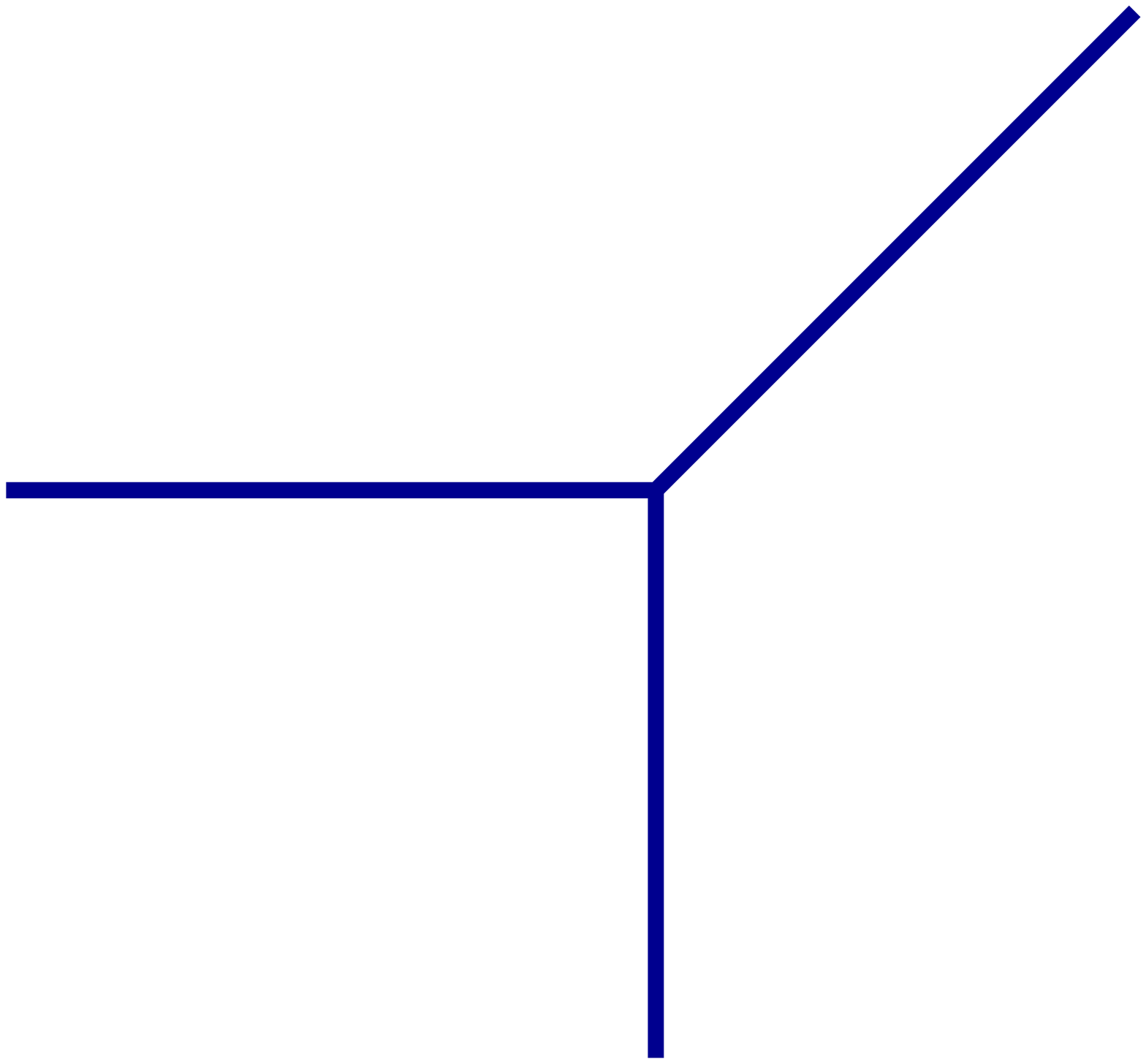}&\hspace*{3cm} &
\includegraphics[width=5cm]{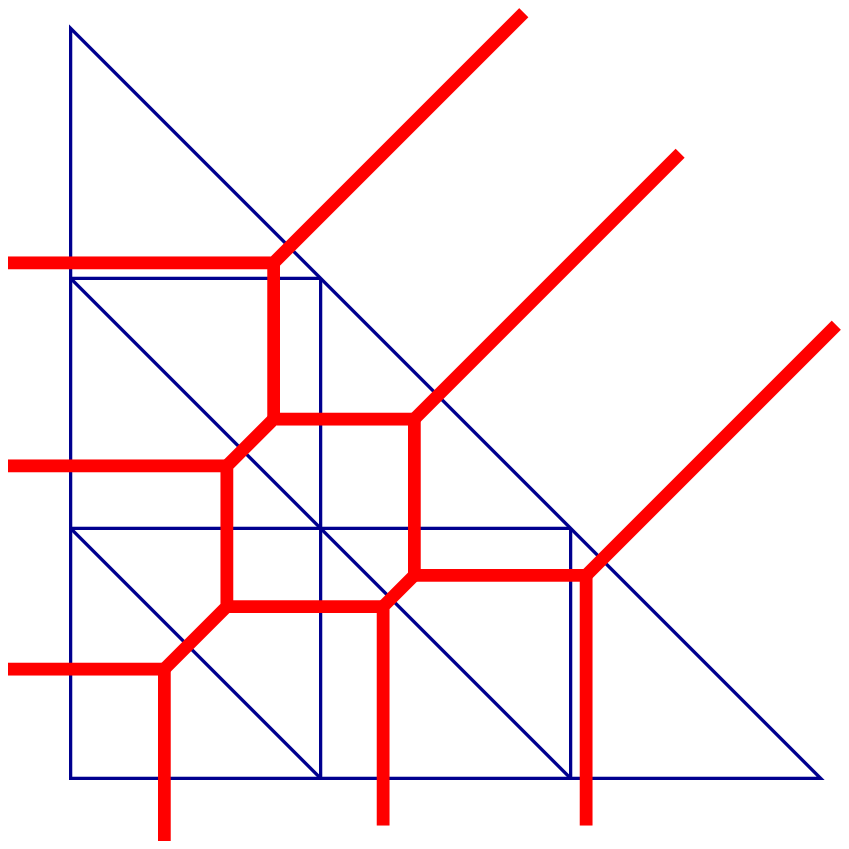}
\end{tabular}
\caption{A tropical line, a tropical cubic with the dual subdivision
  of its Newton polytope.}
\label{F:lineetcubic}
\end{figure}

The linearity domains of $\overline{\nu}$ are the $n$-cells of a convex
polyhedral subdivision $\tau$ of $\Delta$.  The hypersurface
$Z_f^{\trop}$ induces a subdivision $\Xi$ of $\RR^n$. The subdivisions
$\tau$ and $\Xi$ are dual in the following sense (see Figure~\ref{F:lineetcubic}).

There is a one-to-one inclusion reversing correspondence $L$ between
cells of $\Xi$ and cells of $\tau$ such that for any $\xi \in \Xi$,
\begin{enumerate}
\item $\dim L(\xi)=\codim \xi$,
\item the affine supports of $L(\xi)$ and $\xi$ are orthogonal.
\end{enumerate}

We say that $ Z^{\trop}_f $ is {\bf nonsingular} if
its dual subdivision $\tau$ is a primitive triangulation.

\vskip0.6cm

{\bf Complex tropical hypersurfaces.}
In order to prove his correspondence theorem Mikhalkin introduces in
\cite{Mikh05} the complexification of tropical curves.
Let $g(t)=\sum_{r\in R} b_r t^r$ be an element of $\KK^*$. Define the
argument $\arg(g)$ to be the argument of the coefficient of the
smallest power of $t$ appearing in $g$
(i.e. $\arg(g)=\arg(b_{\val(g)}$). Let $\Arg$ be the coordinatewise
argument on $(\KK ^*)^n$. Consider the ``complexification'' of the
valuation 
\[
\begin{array}{rrrl} 
 \VV := V \times \Arg: &(\KK ^*)^n &\longrightarrow &\RR^n
                \times  (S^1)^n\\
                & z  &  \longmapsto &((v(z_1), \dots,
                v(z_n)),(\arg(z_1), \dots, \arg(z_n))).
\end{array}
\]
One can define complex tropical hypersurfaces to be the
closure of the image under $\VV$ of hypersurfaces in  $(\KK ^*)^n$
but it is also convenient to consider the
composition of $\VV$ with the exponential.  Namely, put
\[
\begin{array}{rrrl} 
 \VC : &(\KK ^*)^n &\longrightarrow & (\CC^*)^n\\
                & z  &  \longmapsto &(e^{v(z_1)+i \arg(z_1)},
\dots, e^{v(z_n)+i  \arg(z_n)})
\end{array}
\]

We will call both homeomorphic objects complex tropical hypersurfaces
and use one and the other alternatively depending on the context.

\begin{definition}
  The complex tropical hypersurface $\CC Z_{f,\VC}^{\trop}$ (resp. $\CC
  Z^{\trop}_{f,\VV}$) associated to $f$ is the closure of the image under
  $\VC$ (resp. $\VV$) of the hypersurface $Z_f$:

$$ \CC Z^{\trop}_{f,\VC} := \overline{\VC(Z_f)} \subset (\CC^*)^n,$$

or alternatively

$$ \CC Z^{\trop}_{f,\VV} := \overline{\VV(Z_f)} \subset \RR^n
                \times  (S^1)^n.$$
\end{definition}

\begin{figure}
\begin{tabular}{ccc}
\includegraphics[height=5cm]{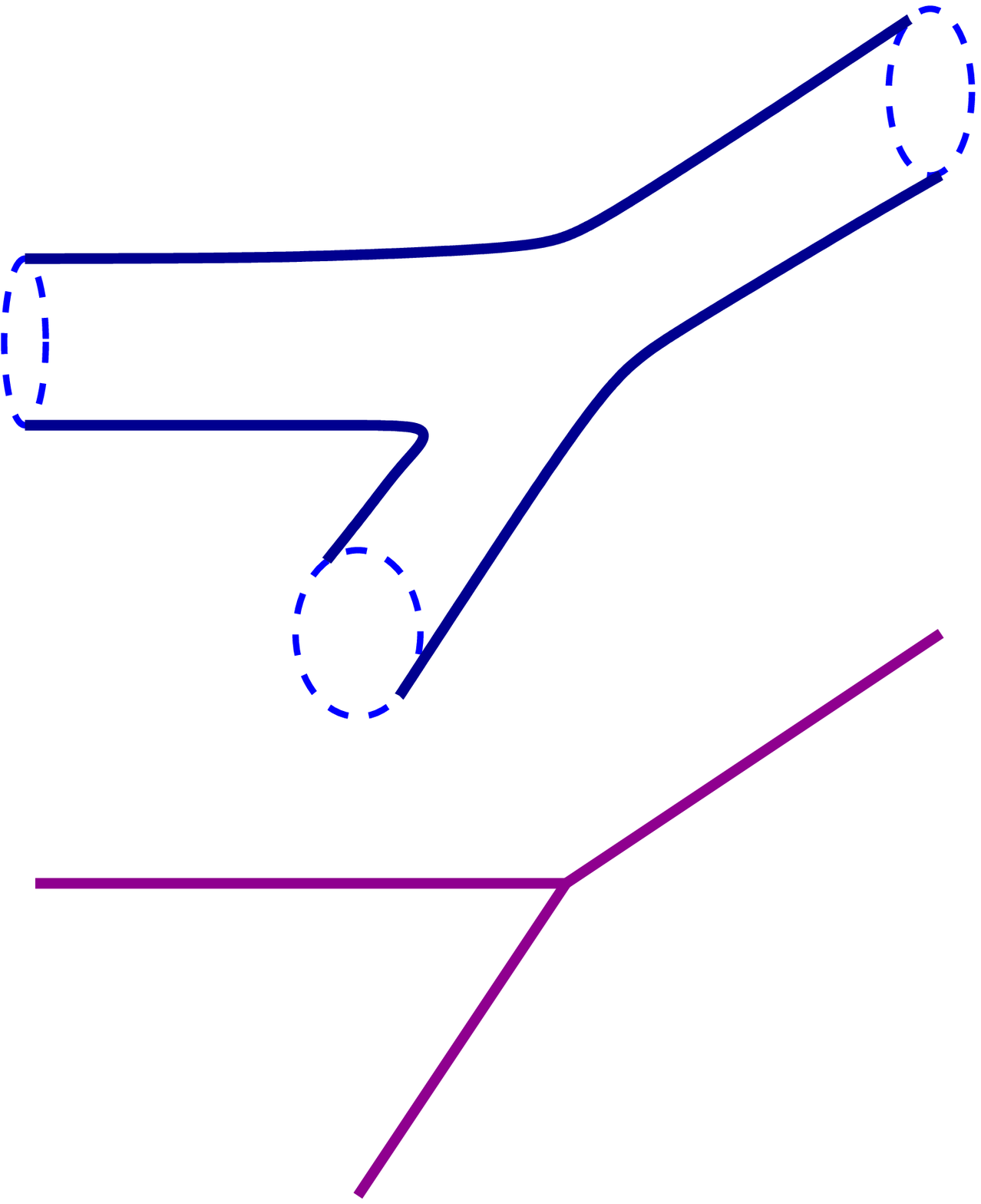}& \hspace*{3cm} &
\includegraphics[height=5cm]{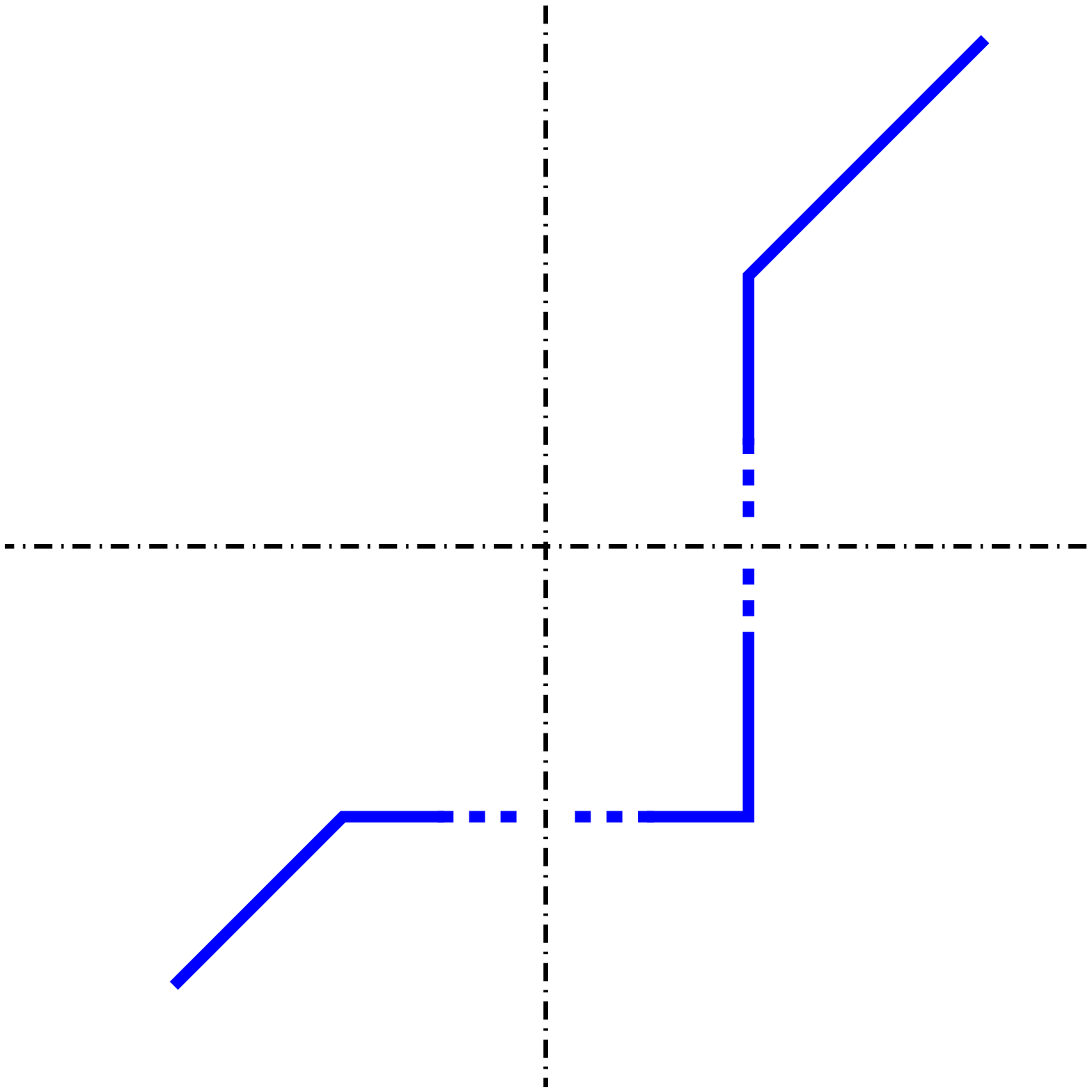}
\end{tabular}
\caption{A complex tropical line, its tropical ``projection'' and a
  real tropical line in the four quadrants.}
\label{F:lines}
\end{figure}

For example, a complex tropical line is homeomorphic to a sphere with
three punctures (see Figure~\ref{F:lines}). 

\vskip0.6cm
{\bf Real tropical Hypersurfaces.} 
Assume from now on that $f=\sum_{\omega \in A} c_\omega z^\omega$ is
real (i.e. all the coefficients $a_r$ of each series $c_\omega =
\sum_{r\in R} a_r t^r $ are real).

\begin{definition}
The real tropical hypersurface $\RR Z_{f,\VC}^{\trop}$ associated to $f$ is 
the intersection in $(\CC^*)^n$ of $\CC Z_{f,\VC}^{\trop}$ with $(\RR^*)^n$:
$$ \RR Z^{\trop}_{f,\VC} := \overline{\VC(Z_f)}\cap(\RR^*)^n,$$
or alternatively
$$ \RR Z^{\trop}_{f,\VV} := \overline{\VV(Z_f)}\cap (\RR^n\times \{0,\pi\}^n).$$
\end{definition}

On Figure~\ref{F:lines} we have pictured the intersection of a real
tropical line with each of the four quadrants.  The sign
$\sign(c_\omega)$ of a coefficient $c_\omega = \sum_{r\in R} a_r t^r$
of $f$ is defined to be the sign of the first coefficient
$a_{\val(c_\omega)}$ of $c_\omega$ (i.e.  $\sign c_\omega$ $:= \sign
a_{\val(c_\omega)}$). From now on we will only consider nonsingular
tropical hypersurfaces. For $\epsilon \in {\{+1,-1\}}^n$, let
$\mathcal{Q}_\epsilon$ be the orthant of $(\RR^*)^n$ which maps to the positive
orthant $(\RR_+^*)^n$ under the isometry $\varphi_\epsilon: \mathcal{Q}_\epsilon
\to (\RR_+^*)^n$ defined by $\varphi_\epsilon(x_1,
\cdots,x_n)=(\epsilon_1 x_1, \cdots,\epsilon_n x_n)$. If
$Z^{\trop}_{f,\VC}$ is nonsingular one can reconstruct $\RR
Z^{\trop}_f \cap \mathcal{Q}_\epsilon$ only from the data of $ Z^{\trop}_f $ and
the collection of signs $\sign c_\omega$ of the coefficients of $f$
(see \cite{Mikhmath.AG/0403015} pp. 25 and 37, \cite{Vir01}, and
\cite{Mik00} Appendix for the case of amoebas).  First one sees that
we can restrict the study to the case of the first orthant $\mathcal{Q}$ since
one can map any $\mathcal{Q}_\epsilon$ to $\mathcal{Q}$ by $\varphi_\epsilon$ and change
$f$ to $f_\epsilon = f(\epsilon_1 x_1, \cdots,\epsilon_n x_n)$.

Let us rather use the map $\VV$ (it suffices to exponentiate
 to switch back to the image under $\VC$). Then $\mathcal{Q}$
corresponds to $\RR^n\times \{(0, \cdots, 0)\}$. Consider the tropical
hypersurface $ Z^{\trop}_f \subset \RR^n$, the induced subdivision
$\Xi_f$ of $\RR^n$ and the dual subdivision $\tau_f$ of its Newton
polytope $\Delta_f$. Let $D_f$ be the sign distribution at the vertices of
$\tau_f$ such that a vertex $\omega$ is labelled with the sign
$\sign(c_\omega)$ of the corresponding coefficient in $f$.

\begin{figure}
\includegraphics[width=4cm]{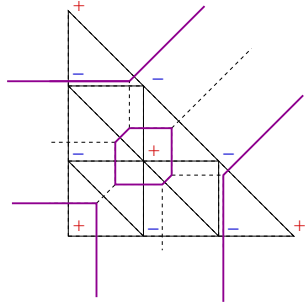}
\caption{A real tropical cubic in the first quadrant, the underlying
  tropical curve (in dotted lines) and its dual
  subdivision.}
\label{F:RT}
\end{figure}

Let $x$ be a point of $Z^{\trop}_f$. The point $(x,(0, \cdots, 0))$ is
in $\RR Z^{\trop}_{f,\VV} \cap (\RR^n\times \{(0, \cdots, 0)\})$ if and only
if $x$ belongs to the closure of an $(n-1)$-cell $\xi$ of $\Xi_f$
which is dual to an edge of $\tau_f$ whose extremities have opposite
signs (see Figure \ref{F:RT} for an example in the case of a curve).

\begin{figure}
\begin{tabular}{cc}
\includegraphics[width=6cm]{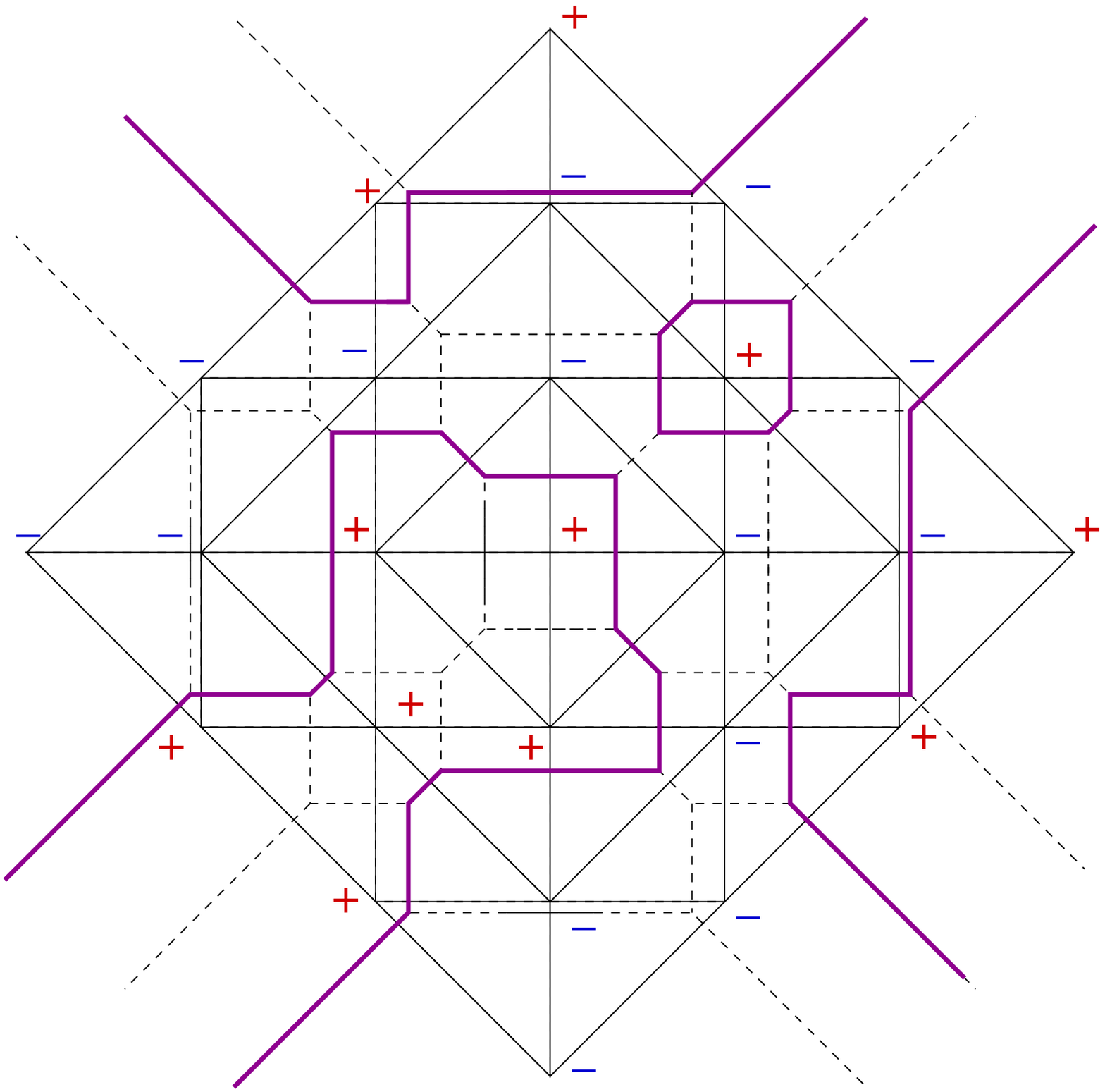}&
\includegraphics[width=6cm]{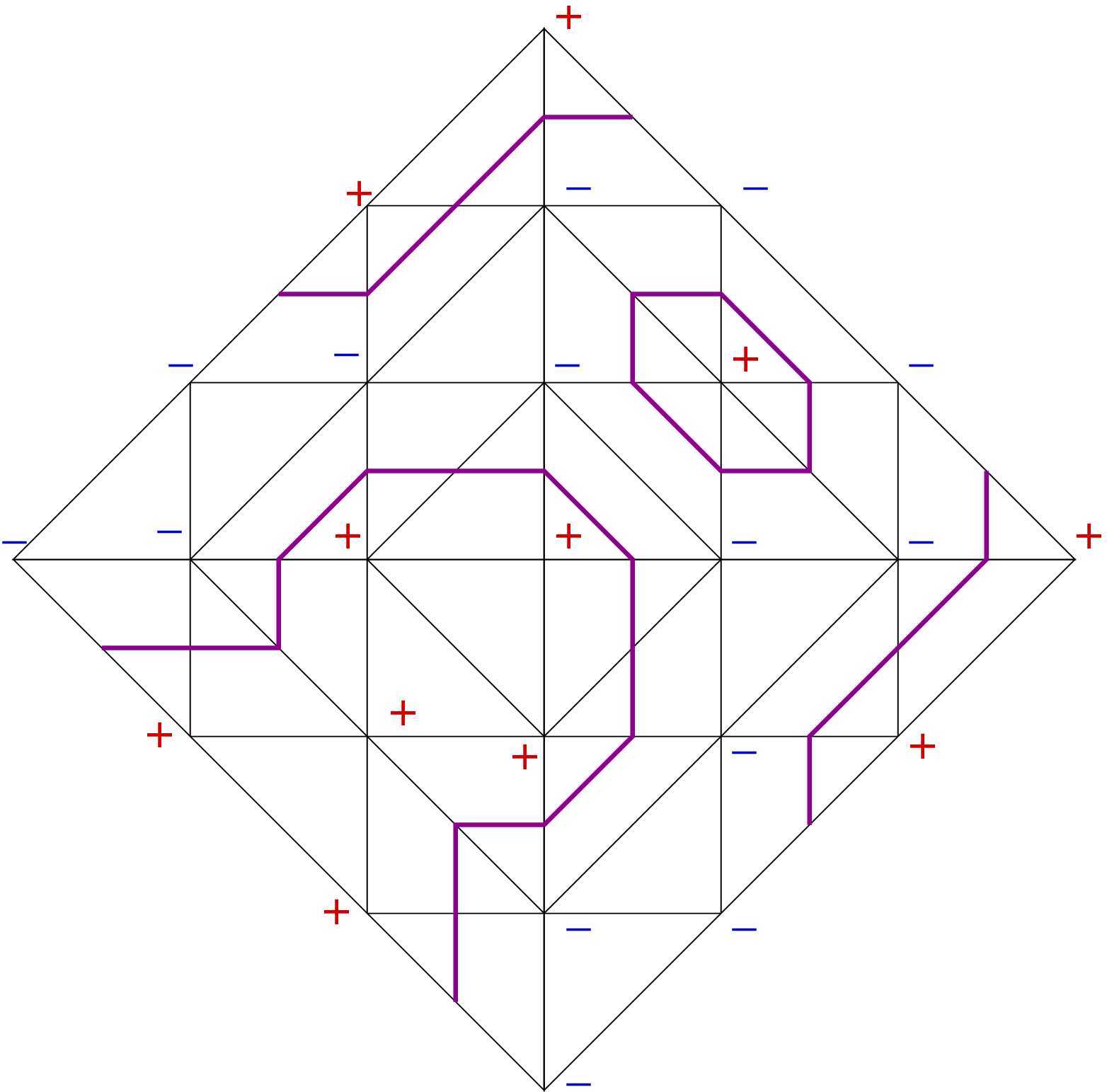}
\end{tabular}
\caption{An attempt at representing a real nonsingular tropical
  cubic with the underlying symmetrized tropical curve (in dotted
  lines) and corresponding $T$-curve.}
\label{F:patchtrop}
\end{figure}

Then one sees that one can associate a unique primitive
$T$-hypersurface to a nonsingular real tropical hypersurface. Namely
take the $T$-hypersurface $H$ constructed from $\Delta_f, \tau_f$ and
$D_f$. Clearly, there exists a homeomorphism $h:(\RR^*)^n \to (\Delta \setminus
\partial\Delta) ^*$ such that $h(\RR Z^{\trop}_f)= H \cap (\Delta \setminus
\partial\Delta)^*$. On Figure~\ref{F:patchtrop} we tried to illustrate
this homeomorphism in the case of curves.

The results in this paper are nicer and much easier to state when one
compactifies the tropical hypersurfaces. Let us describe the natural
compactification of a tropical hypersurface. Recall that $\CC
Z^{\trop}_{f,\VC} = \overline{\VC(Z_f)}$ is a subset of the torus
$(\CC^*)^n$.  Consider the usual compactification of $(\CC^*)^n$ into
the toric variety $X_{\Delta_f}$ associated to the Newton polytope
$\Delta_f$ of $f$.  Let $\iota: (\CC^*)^n \hookrightarrow
X_{\Delta_f}$ be the corresponding inclusion. 
 We define the compactification $\overline{\CC
  Z^{\trop}_{f}}$ to be the closure of $\iota(\CC
Z^{\trop}_{f,\VC})$ in $X_{\Delta_f}$. 
 Note that the stratification of $X_{\Delta_f}$ into orbits of
the action of $(\CC^*)^n$ defines a natural stratification of
$\overline{\CC Z^{\trop}_{f}}$.

We sum up natural maps in  the following commutative diagram.
\[
  \xymatrix{
  \RR Z^{\trop}_{f,\VV}  \ar[r]^\sim\ar@{^{(}->}[d] & \RR Z^{\trop}_{f,\VC}
   \ar@{^{(}->}[r]\ar@{^{(}->}[d] & {(\RR^{\relax *})^n}\ar@{^{(}->}[d] \ar@{^{(}->}[r]^{\iota_\RR} & \RR X_{\Delta_f} \ar@{^{(}->}[d] \\
   {\CC Z^{\trop}_{f,\VV}} \ar[r]^\sim & {\CC Z^{\trop}_{f,\VC}} \ar@{^{(}->}[r] & (\CC^{\relax *})^n\ar@{^{(}->}[r]^\iota & X_{\Delta_f}
  }
\]
Define $\overline{\RR Z^{\trop}_f}$ to be the
intersection of $\overline{{\CC Z^{\trop}_{f}}}$ with the
real part $\RR X_{\Delta_f}$ of $X_{\Delta_f}$. Clearly $\overline{\RR
  Z^{\trop}_f}$ is also the closure of $\iota_\RR(\RR
Z^{\trop}_f)$ in $\RR X_{\Delta_f}$. One can see that the natural
stratification of $\overline{\RR Z^{\trop}_f}$ induced
by the torus action corresponds to the stratification of the
$T$-hypersurface $H$ induced by the face complex of $\Delta_f$ which
will be used in the proof of Theorem~\ref{theprop}.

\subsection{Tropical Statements}\label{S:tropstate}
We can now state our results in the tropical language. 
Let  $\Delta$ be a polytope corresponding to a nonsingular toric variety.

\begin{bbthm}\label{Th:sigmachitrop}
  Let $\overline{\RR Z^{\trop}_f}$ be a compactified
  nonsingular real tropical hypersurface with Newton polytope
  $\Delta_f=\Delta$ then

$$\chi(\overline{\RR Z^{\trop}_f}) = \sigma(Z),$$

where $Z$ is a 
smooth complex algebraic hypersurface in
$X_{\Delta_f}$ with Newton polytope $\Delta_f$ and $\sigma(Z)=
\sum_{p+q=0\> [2]} {(-1)}^p h^{p,q}(Z)$.
\end{bbthm}

Theorem~\ref{Th:sigmachitrop} is equivalent to Theorem~\ref{theprop}
about Euler characteristic of primitive $T$-hypersurfaces and is one
the of main results of this paper.

One can define a maximal real tropical hypersurface in the following
way. Let $f$ be a polynomial over $\KK$ such that every coefficient of
$f$ is a series with only real coefficients. Let $\Delta_f$ be the
Newton polytope of $f$. The tropical hypersurface $\overline{\RR
  Z^{\trop}_f}$ is {\bf maximal} if 
\[
b_*(\overline{\RR
  Z^{\trop}_f};\ZZ_2) = b_*(\overline{\CC
  Z^{\trop}_f};\ZZ_2). 
\]
 
We can then prove existence results for maximal real tropical
surfaces.


\begin{bbthm}\label{Th:P2trop}
Let $\Delta$ be a $3$-dimensional Nakajima polytope corresponding to a
nonsingular toric variety $X_\Delta$.
Then, there exists a maximal real tropical surface
in $X_\Delta$ with the Newton polytope $\Delta$.
\end{bbthm}

Theorem\ref{Th:P2trop} will be proved in Section~\ref{Mhypex} (See Theorem~\ref{P2})
as well as the following corollary (See Corollary~\ref{hirz}).

For a nonnegative integer $\alpha$ and positive integers $m$ and $n$
denote by $\delta_\alpha^{m,n}$ the polygon having the vertices
$(0,0)$, $(n + m \alpha,0)$, $(0,m)$, and $(n,m)$ in $\RR^2$.  The
toric variety associated with $\delta_\alpha^{m,n}$ is a rational
ruled surface $\Sigma_\alpha$.  Consider now the truncated cylinder
$\Delta^{l,m,n}_\alpha$ with base $\delta_\alpha^{m,n}$ whose vertices are
$$\displaylines{(0,0,0),(n + m \alpha,0,0),(0,m,0),(n,m,0), \cr
  (0,0,l), (n + m \alpha,0,l),(0,m,l), \; \text{\rm and} \;
  (n,m,l),}$$ 
where $l$ is a positive integer.  The toric variety
  $X_{\Delta^{l,m,n}_\alpha}$ associated with $\Delta^{l,m,n}_\alpha$
  is $\Sigma_\alpha \times \CC P^1$. Since $\Delta^{l,m,n}_\alpha$ is,
  up to exchanging two coordinates, a Nakajima polytope, we have the
  foolowing corollary (See Corollary~\ref{hirz}).

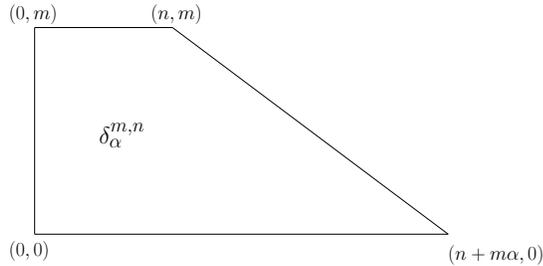
\begin{figure}[htbp]
\begin{center}
\resizebox {6 cm}{!}{
\input{dessins/delta.pstex_t}
}
\caption{Polygon $\delta_\alpha^{m,n}$.}
\label{figure:delta}
\end{center}
\end{figure}

\begin{large}
\begin{figure}[htbp]
\begin{center}
\resizebox {6 cm}{!}{
\input{dessins/P.pstex_t}
}
\caption{Polytope $\Delta^{l,m,n}_\alpha$.}
\label{figure:P}
\end{center}
\end{figure}
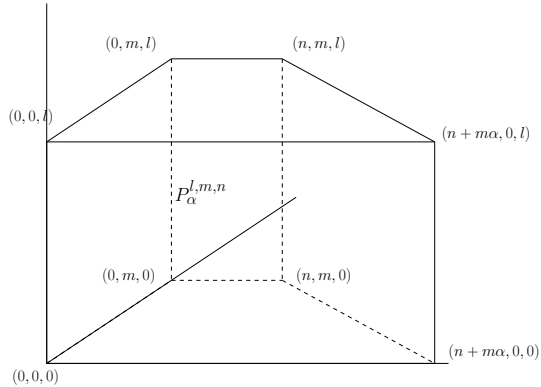
\end{large}

\begin{bbcor}\label{Th:hirztrop}
  For any nonnegative integer $\alpha$ and any positive integers $m$,
  $n$ and $l$, there exists a (real) maximal nonsingular tropical
  surface $\overline{\RR Z^{\trop}_{f}}$ with the Newton polytope
  $\Delta^{l,m,n}_\alpha$ in $\RR\Sigma_\alpha \times \RR P^1$.
\end{bbcor}

\subsection{Some facts about triangulations of lattice polytopes}

\subsubsection{Ehrhart Polynomial}\label{ehrpol}

In order to compute 
both sides of equality~(\ref{sc}),
 namely the Euler characteristic of the real part of a primitive
$T$-hypersurface and the signature of its complex part, we use some
combinatorics of triangulations of lattice polytopes. We first recall
the results of E. Ehrhart on the number of integer points in a lattice
polytope.  Ehrhart showed that the number of integer points in the
multiple $\lambda \cdot \Delta$ of a (not necessarily convex) polytope
$\Delta$, where $\lambda$ is a positive integer, is a polynomial in
$\lambda$ (see \cite{Ehr} or \cite{Ehr67}).  We denote by $l(\Delta)$
and $l^*(\Delta)$ the numbers of integer points in $\Delta$ and in the
interior of $\Delta$, respectively.

\begin{bbthm}[E. Ehrhart]
  Let $\Delta$ be a polytope with integer vertices.  Then, the numbers
  $l(\lambda \cdot \Delta)$ and $l^*(\lambda \cdot \Delta)$ are
  polynomials in $\lambda$ of degree $\dim \Delta$.  Denote them
  respectively by $Ehr_\Delta(\lambda)$ and $Ehr^*_\Delta(\lambda)$.
  They satisfy the reciprocity law
$$ (-1)^{dim \Delta} Ehr_\Delta(-\lambda)=Ehr^*_\Delta(\lambda) .$$
\end{bbthm}

One often considers the Ehrhart series
$$SE_p(t)=\sum_{\lambda=0}^{\infty} Ehr_\Delta(\lambda) t^{\lambda}.$$
Put $Q_\Delta(t) =(1-t)^{d+1} SE_p(t)$, where $d$ is the dimension of
$\Delta$. 
 We define the numbers $\Psi_j$ to
be the coefficients of $Q_\Delta(t)$ :

$$Q_\Delta(t)= \sum_{j=0}^\infty \Psi_j t^j .$$

 In fact, $Q_\Delta(t)$ is a polynomial of degree $d$ (see
\cite{bri} or lemma \ref{psi}).
Let $a_i^{\Delta}$ be the coefficient of $\lambda^i$ in  $Ehr_\Delta(\lambda)$.
The following formula can be found in Section~4.1 of \cite{DCZ}  (See also the Appendix).

\begin{bblem}\label{psi}
  One has $$\Psi_j= \sum_{i=0}^d (\sum_{n=0}^j
  (-1)^{j-n}\bin{d+1}{j-n} n^i) a_i^{\Delta},$$
  and $\Psi_j= 0$ for $ j \ge
  d+1$.
\end{bblem}

\subsubsection{Number  of  simplices  in  a  primitive
  triangulation}\label{nbs}

Another usefull formula 
expresses the numbers of simplices of any dimension in the primitive
triangulation $\tau$ of $\Delta$ in terms of the coefficients of the
Ehrhart polynomial of $\Delta$.  In the case of a primitive
triangulation, these numbers happen not to depend on the primitive
triangulation chosen.  This statement can be found in \cite{Dai} and
in 
 Appendix~\ref{ksimplices} 
 we provide a proof in the case of convex
triangulations.

Let $nbs_{r}^{F}$ be the number of
$r$-dimensional simplices in a primitive
triangulation of $F$ which are contained in
the interior of $F$, and let $i$ be the dimension
of $F$.
Let $S_2$ be the Stirling number of the second kind defined by
$S_2(i,j)=1/(j)! \sum_{k=0}^{j} (-1)^{j-k}\bin{j}{k} k^i$.

\begin{bbprop}\label{nbsrformula}
Under the above hypotheses we have,
$$N_{r}^{F}=\sum_{l=r+1}^{i+1}
(-1)^{i-l+1}\cdot r! S_2(l,r+1) \cdot a^{F}_{l-1}.$$
\end{bbprop}

\subsection{Danilov and Khovanskii formulae}\label{S:DK}

V. Danilov and A. Khovanskii \cite{DanKho} computed the Hodge numbers
of a smooth hypersurface in a toric variety $X_\Delta$ in terms of the
polytope $\Delta$.
Recall that,
 for a face $F$ of $\Delta$, the
coefficient of the term of degree $l$ of the Ehrhart polynomial of
$F$ is denoted by $a_l^F$.
  Let $\mathcal{F}_i(\Delta)$ be the set of $i$-dimensional faces
of $\Delta$, and $f_i(\Delta)$ be the cardinality of
$\mathcal{F}_i(\Delta)$.  A $d$-dimensional {\bi simple} polytope is
one for which each vertex is contained in exactly $d$ edges.

\begin{bbthm}[V. Danilov, A. Khovanskii]\label{DK}
  Let $\Delta$ be a simple polytope of dimension $d$, and $Z$ be a
  smooth or quasi-smooth\footnote{See \cite{Dan78} and
  Remark~\ref{quasismooth}} algebraic hypersurface with Newton
  polytope $\Delta$ in $X_\Delta$.  Then, for $p\neq \frac{d-1}{2}$
\begin{eqnarray*}
h^{p,p}(Z) & = & {(-1)}^{p+1} \sum_{i=p+1}^{d}{(-1)}^i \bin{i}{p+1}
f_i(\Delta)\\
 h^{\frac{d-1}{2},\frac{d-1}{2}}(Z) & = & {(-1)}^{\frac{d+1}{2}}
\sum_{i=\frac{d+1}{2}}^{d}{(-1)}^i \bin{i}{\frac{d+1}{2}} f_i(\Delta)
- \sum_{i=\frac{d+1}{2}}^{d} \sum_{F \in \mathcal{F}_i(\Delta)} {(-1)}^i
\Psi_{\frac{d+1}{2}}(F)\\
h^{p,d-1-p}(Z) & = & (-1)^{d} \sum_{i=p+1}^{d} \sum_{F
\in \mathcal{F}_i(\Delta)}{(-1)}^i \Psi_{p+1}(F)\\
h^{p,q} & = & 0 \mbox{ if } q \neq p \; \mbox{or} \; p \neq d-1-p,
\end{eqnarray*}
where 
$\Psi_{p+1}(F)=\sum_{\alpha=1}^{i+1}\sum_{a=0}^{p+1} (-1)^a
 \bin{i+1}{a}(p+1-a)^{\alpha -1} a^{F}_{\alpha -1 } $.
\end{bbthm}

\section{Primitive $T$-hypersurfaces}\label{sigmachi}

\subsection{Statement}

Let $\Delta$ be a $d$-dimensional polytope corresponding to a
nonsingular toric variety.  Let $Z$ be a primitive $T$-hypersurface with
Newton polytope $\Delta$ in $X_\Delta$, and $\tau$ be
 a primitive triangulation 
from which one can construct $Z$ via the combinatorial Viro method.
 Notice that since
$X_\Delta$ is smooth the $T$-hypersurface $Z$ is also smooth.
Define $\sigma(Z)$ by 

$$\sigma(Z):= \sum_{p+q=0\> [2]} {(-1)}^p h^{p,q}(Z).$$

If $d$ is odd and $Z$ is smooth, the intersection form 
$$\iota: H_{d-1}(Z;\ZZ)/tors
\times H_{d-1}(Z;\ZZ)/tors \rightarrow \ZZ$$ 
is a symmetric bilinear
form and its signature is called the signature of $Z$ and is equal to
$\sigma(Z)$. (See, for example, \cite{GH} Chapter~0, Section~7 p. 126).
We have the following statement.

\begin{bbthm}\label{theprop}
If $Z$ is a primitive real algebraic $T$-hypersurface in a
 nonsingular toric variety $X_\Delta$,  then
$$\chi(\RR Z) = \sigma(Z),$$
where $\chi(\RR Z)$ is the Euler characteristic of the real part $\RR
Z$ of $Z$.
\end{bbthm}

\begin{rem}\label{quasismooth}
For simplicity we state Theorem~\ref{theprop}  for a
nonsingular ambiant toric variety.  In fact, we only need to require
that the polytope $\Delta$ be simple (see Subsection~\ref{S:DK}).  In
this case $X_\Delta$ is quasi-smooth (i.e. locally isomorphic to a
toric variety defined by a simplicial cone; see
\cite{Dan78}~Section~14). Then $Z$ is also quasi-smooth (see
\cite{Dan78}~Section~13
), the Hodge numbers can be calculated by Danilov and Khovanskii formulae (see Theorem~\ref{DK}) and the proof is verbatim.  
\end{rem}

\begin{rem}
The signature $\sigma(Z)$ is the same for any generic hypersurface $Z$
 with Newton polytope $\Delta$ in $X_\Delta$ but if
 $Z$ is not a primitive $T$-hypersurface $\chi(\RR Z)$  needs not
 be equal to the signature (and very often is not).  
\end{rem}

Recall that $l^*(\Delta)$ is the number of integer points in the
interior of $\Delta$. Here and further on we write $b_i(X)$ for
$b_i(X;\ZZ_2)$.

\begin{bbcor}\label{max}
Assume that $\Delta$ is
a $3$-dimensional simple polytope and $Z$ is
a primitive T-surface in $X_\Delta$ with Newton polytope $\Delta$.
If the number
$b_0(\RR Z)$ of connected components of $\RR Z$ is at least
$l^*(\Delta) + 1$, then
$b_0(\RR Z)=l^*(\Delta) + 1$ and $Z$ is maximal.
\end{bbcor}

\proof First, note that $l^*(\Delta)=h^{2,0}(Z)$.
Then $\chi(\RR Z) \ge 2h^{2,0}(Z) + 2 - b_1(\RR Z)$.
Now using
the equalities $\chi(\RR Z) =
\sigma( Z)=2h^{2,0}(Z) + 2 - h^{1,1}(Z)$
one gets $b_1(\RR Z) \ge h^{1,1}(Z)$. Furthermore $h^{1,0}(Z)=h^{0,1}(Z)=0$
(see Theorem \ref{DK}), and thus $b_*(Z; \CC ) =2h^{2,0}(Z) + 2 +
h^{1,1}(Z)$.
Hence $b_*(\RR Z; \ZZ_2) \geq b_*(Z;  \CC) = b_*(Z; \ZZ_2)$.
The Smith-Thom inequality implies that
$b_*(\RR Z;\ZZ _2) =  b_*(Z;\ZZ_2)$,
and thus,
that $b_1(\RR Z)=h^{1,1}(Z)$ and $b_0(\RR Z)=h^{2,0}(Z)+1$.\CQFD

\subsection{ Proof of Theorem \ref{theprop}}

If $d$ is even,
Theorem~\ref{theprop} is
straightforward. Indeed, in this case $Z$ is a (smooth) odd dimensional
hypersurface, so $\chi(\RR Z)=0$.
On the other hand, we have the
equality $h^{p,q}(Z)=h^{d-1-p,d-1-q}(Z)$
for any $p$ and $q$,
and $d-1$ is odd.
Thus, $\sigma(Z)= \sum_{p+q=0\> [2]} {(-1)}^p h^{p,q}(Z) = 0$.

Assume now that $d$ is odd.
\begin{bblem}\label{exprsigchi}
We have
\begin{eqnarray*}
 \chi(\RR Z) = \sum_{i=1}^{d}\sum_{F\in
\mathcal{F}_i(\Delta)} \sum_{l=2}^{i+1} \chi_{l,i+1}
a_{l-1}^F & \mbox{and} &
\sigma(Z) = \sum_{i=1}^{d}\sum_{F\in
\mathcal{F}_i(\Delta)} \sum_{l=2}^{i+1} \sigma_{l,i+1}
a_{l-1}^F,\\
\end{eqnarray*}
where 
\begin{eqnarray*}
\chi_{l,i+1}  & = & (-1)^{i-l+1} \sum_{k=1}^{i}
 \frac{(2^i-2^{i-k})}{k+1} \sum_{m=0}^{k+1}
(-1)^{m}\bin{k+1}{m} m^l,\\
\sigma_{l,i+1} & = &
 \sum_{p=0}^{d-1} {(-1)}^i
  \sum_{q=0}^{p+1} (-1)^{p+1-q}
\bin{i+1}{q}(p+1-q)^{l -1}.
\end{eqnarray*}
\end{bblem}
\proof The triangulation $\tau^*$ induces a cell decomposition
$\mathcal{D}$ of $\widetilde{H}$. Let $\widetilde{I_F}$ be the image
in $\widetilde{\Delta}$ of the union of the symmetric copies of the
interior of a face $F$ and $\mathcal{D}(\widetilde{I_F})$ be the set
of cells of $\mathcal{D}$ contained in $\widetilde{I_F}$.  Put $\chi_F
= \sum_{\delta \in \mathcal{D}(\widetilde{I_F}) } (-1)^{dim(\delta)}$.
Then $\chi(\RR Z)=\sum_{i=1}^{d}\sum_{F\in \mathcal{F}_i(\Delta)}
\chi_F$.  According to Proposition \ref{itenprop}, if $Q$ is a
$k$-simplex of $\tau$ contained in $j$ coordinate hyperplanes then the
union of the symmetric copies of $Q$ contains exactly
$(2^d-2^{d-k})/2^j$ cells of dimension $k-1$. An $i$-face $F$ of
$\Delta^*$ contained in $j$ coordinate hyperplanes is identified with
$2^{d-i-j}-1$ other copies of $F$ when passing from $\Delta^*$ to
$\widetilde{\Delta}$. Thus, the number of $(k-1)$-cells in
$\mathcal{D}(\widetilde{I_F})$ is equal to
$\frac{(2^d-2^{d-k})}{2^{d-i}} N^F_k$, where $N_k^F$ is the number
of $k$-simplices in the interior of $F$. Thus,
\begin{equation}
 \chi_F  =  \sum_{k=1}^{\dim F} (-1)^{k-1} (2^{\dim F} -2^{\dim F-k}) N_{k}^F
 \; .  
\end{equation}
According to
Proposition \ref{nbsrformula} of
Section \ref{nbs}
$$N_{k}^{F}=\sum_{l=k+1}^{\dim F+1}
k!S_2(l,k+1)(-1)^{\dim F-l+1}a^{F}_{l-1},$$
\noindent
where
$S_2(i,j)=1/j! \sum_{k=0}^{j}(-1)^{j-k}\bin{j}{k} k^i$ is the
Stirling number of the second kind. Then one has
\begin{eqnarray*}
\chi_F & = & \sum_{k=1}^{\dim F} (-1)^{k-1} (2^{\dim F} -2^{\dim F-k}) N_{k}^F
 \\
& = & \sum_{k=1}^{\dim F} (-1)^{\dim F -1} \frac{2^{\dim F} -2^{\dim F-k}}{k+1}
\sum_{l=k+1}^{\dim F +1} (-1)^l a_{l-1}^F \sum_{m=0}^{k+1}(-1)^m
 \bin{k+1}{m} m^l \; .  
\end{eqnarray*}

By Lemma~\ref{vanlint}, Equality~(\ref{nbsurj})  (see Appendix), the sum on $l$ can be
taken from $l=2$ instead of $l=k+1$ and one gets
\begin{equation*}
\chi_F  =  \sum_{l=2}^{\dim F +1} (-1)^l a_{l-1}^F \sum_{k=1}^{\dim F}
(-1)^{\dim F-1} \frac{2^{\dim F} -2^{\dim F -k}}{k+1} \sum_{m=0}^{k+1}(-1)^m
 \bin{k+1}{m} m^l \; ,
\end{equation*}
which finishes the proof of the formula for $\chi(\RR Z)$.

To  compute $\sigma(Z)$, we use the Danilov and  Khovanskii formulae
(see Theorem \ref{DK}).
One obtains the following expression
for $\sigma (Z)$:

\begin{small}
\begin{eqnarray*}
\sigma(Z) & = &  \sum_{p=0}^{d-1} \sum_{i=p+1}^{d}{(-1)}^i \sum_{F\in
\mathcal{F}_i(\Delta)}
\Bigg(- \bin{i}{p+1} + {(-1)}^{p+1} \sum_{l=1}^{i+1}\sum_{q=0}^{p+1} (-1)^q
\bin{i+1}{q}(p+1-q)^{l -1} a^{F}_{l -1} \Bigg).
\end{eqnarray*}
\end{small}

Consider $\sigma(Z)$ as an affine polynomial in the variables $a^{F}_{l}$.
Denote by $\sigma^{cst}$ the constant term of $\sigma(Z)$,
and by $\sigma^1$ the sum of monomials
in variables $a^{F}_{0}$.
We have
$$\sigma^{cst}= - \sum_{p=0}^{d-1} \sum_{i=p+1}^{d}{(-1)}^i \bin{i}{p+1} f_i(\Delta) $$
and
$$\sigma^1 = \sum_{p=0}^{d-1} \sum_{i=p+1}^{d}{(-1)}^i {(-1)}^{p+1}
\sum_{F\in\mathcal{F}_i(\Delta)} \sum_{q=0}^{p+1} (-1)^q \bin{i+1}{q}
a^F_0 .$$

\begin{bblem}
One has
$\sigma^{cst}+\sigma^1=0$.
\end{bblem}

\proof
Since $ a^0_F= 1$, we obtain
\begin{eqnarray*}
\sigma^1 &= &\sum_{p=0}^{d-1} \sum_{i=p+1}^{d}{(-1)}^i {(-1)}^{p+1}
\sum_{F\in\mathcal{F}_i(\Delta)} \sum_{q=0}^{p+1} (-1)^q \bin{i+1}{q} \\
& = & \sum_{p=0}^{d-1} \sum_{i=p+1}^{d}{(-1)}^i {(-1)}^{p+1} f_i(\Delta)
\sum_{q=0}^{p+1} (-1)^q \bin{i+1}{q}.
\end{eqnarray*}
Furthermore,
$\sum_{q=0}^{p+1} (-1)^q \bin{i+1}{q}  = (-1)^{p+1}\bin{i}{p+1}$,
and thus,
\begin{eqnarray*}
\sigma^1 & = &\sum_{p=0}^{d-1}
\sum_{i=p+1}^{d}{(-1)}^i \bin{i}{p+1} f_i(\Delta).  
\end{eqnarray*}
\proofend

For any $p \ge i$, $\sum_{q=0}^{p+1} (-1)^q
\bin{i+1}{q}(p+1-q)^{l -1} =\sum_{q=0}^{i+1} (-1)^q
\bin{i+1}{q}(p+1-q)^{l -1}$  which is zero by Lemma~\ref{vanlint}
 since $l\le i+1$. Thus one gets 
\begin{equation}\label{sli}
\sigma_{l,i+1}  = \sum_{p=0}^{i-1} {(-1)}^i {(-1)}^{p+1}
\sum_{q=0}^{p+1} (-1)^q
\bin{i+1}{q}(p+1-q)^{l -1}.
\end{equation}
\proofend

We can then prove the following lemma.

\begin{bblem}\label{sigchirec}
For $2 \le l \le i+1$, one has the following equalities
\begin{eqnarray*}
\sigma_{l,i+2}& = &-2\sigma_{l,i+1},\\
 \chi_{l,i+2}  & = & -2 \chi_{l,i+1}.
\end{eqnarray*}
\end{bblem}

\proof
To prove the first equality,
write  $\bin{i+2}{q}=\bin{i+1}{q-1} + \bin{i+1}{q}$
to get
$$\displaylines{
\sigma_{l,i+2} = {(-1)}^{i+1}\sum_{p=0}^{i}
 {(-1)}^{p+1} \sum_{q=1}^{p+1} (-1)^q \bin{i+1}{q-1}(p+1-q)^{l
 -1}\cr
+ {(-1)}^{i+1}\sum_{p=0}^{i}
 {(-1)}^{p+1} \sum_{q=0}^{p+1} (-1)^q \bin{i+1}{q}(p+1-q)^{l
 -1}.}$$

Notice that $\sum_{q=1}^{1} (-1)^q \bin{i+1}{q-1}(p+1-q)^{l -1}$ is
$0$.  Similarly, in the second term,\\ $\sum_{q=0}^{i+1} (-1)^q
\bin{i+1}{q}(i+1-q)^{l -1}$ makes no contribution,
 by Lemma \ref{vanlint}.

Then,
$$\displaylines{
\sigma_{l,i+2} = {(-1)}^{i+1}\sum_{p=1}^{i}
 {(-1)}^{p+1} \sum_{q=1}^{p+1} (-1)^q \bin{i+1}{q-1}(p+1-q)^{l
 -1}\cr
+ {(-1)}^{i+1}\sum_{p=0}^{i-1}
 {(-1)}^{p+1} \sum_{q=0}^{p+1} (-1)^q \bin{i+1}{q}(p+1-q)^{l
 -1}.}$$

So, with the changes of indices $ c=q-1$ and $ r=p-1$ one gets
\begin{eqnarray*}
\sigma_{l,i+2}&=&  {(-1)}^{i+1}\sum_{r=0}^{i-1}
 {(-1)}^{r} \sum_{c=0}^{r+1} (-1)^{c+1} \bin{i+1}{c}(r+1-c)^{l
 -1} \\
& & + {(-1)}^{i+1}\sum_{p=0}^{i-1}
 {(-1)}^{p+1} \sum_{q=0}^{p+1} (-1)^q \bin{i+1}{q}(p+1-q)^{l
 -1}\\
&=& -2 \sigma_{l,i+1}.
\end{eqnarray*}

To prove the second equality stated in the lemma,
use Lemma \ref{vanlint} to
 get rid of the term with $k=i+1$ in the expression of $\chi_{l,i+2}\,$:

\begin{eqnarray*}
\chi_{l,i+2}  & = & (-1)^{i-l+2}  \sum_{k=1}^{i+1}
 \frac{(2^{i+1}-2^{i+1-k})}{k+1} \sum_{m=0}^{k+1}
(-1)^{m}\bin{k+1}{m} m^l \\
 & = &  2 (-1)^{i-l}  \sum_{k=1}^{i}
 \frac{(2^{i}-2^{i-k})}{k+1} \sum_{m=0}^{k+1}
(-1)^{m}\bin{k+1}{m} m^l \,,
\end{eqnarray*}

which finishes the proof.
\CQFD

\begin{bblem}\label{sigeqchi}
For $2 \le l \le i+1 \le d+1$, we have $\chi_{l,i+1}=\sigma_{l,i+1}$.
\end{bblem}

\proof  By Lemma~\ref{sigchirec} it is enough to
prove the equality
$\chi_{l,l}=\sigma_{l,l}$.

By Formula~(\ref{sli}) one has, 
\begin{eqnarray}\label{sll}
\sigma_{l,l} & = &
 \sum_{p=0}^{l-2} {(-1)}^{l-1}
 \sum_{q=0}^{p+1} (-1)^{p+1-q}
\bin{l}{q}(p+1-q)^{l -1}\label{sll1}\\
& = &\sum_{b=0}^{l} {(-1)}^{l-1}
 \sum_{q=0}^{l-b} (-1)^b
\bin{l}{q}(b)^{l-1}\label{sll2}\\
& = &  {(-1)}^{l-1} \sum_{b=1}^{l}\sum_{j=b}^{l}
 (-1)^{b} \bin{l}{j}b^{l -1},\label{sll3}
\end{eqnarray} 

where we pass from (\ref{sll1}) to  (\ref{sll2})
 by setting $b=p+1-q$ and from (\ref{sll2}) to (\ref{sll3}) by putting $j=l-q$.

One sees that
\begin{eqnarray*}
\chi_{l,l} &= & \sum_{k=1}^{l-1}
 \frac{2^{l-1}-2^{l-1-k}}{k+1} \sum_{m=0}^{k+1}
(-1)^{m}\bin{k+1}{m} m^l \\
&= &\sum_{k=0}^{l-1}
 \frac{2^{l-1}-2^{l-1-k}}{k+1} \sum_{m=1}^{k+1}
(-1)^{m}\bin{k+1}{m} m^l \,, 
\end{eqnarray*}
just noticing that the change of range does not affect the sum.
Since $\bin{k+1}{m}\frac{m^l}{k+1}= \bin{k}{m-1} m^{l-1}$,

\begin{eqnarray*}
\chi_{l,l}&= &\sum_{k=0}^{l-1}
 (2^{l-1}-2^{l-1-k}) \sum_{m=1}^{k+1}
(-1)^{m}\bin{k}{m-1} m^{l-1} \\
& = &  2^{l-1} \sum_{k=0}^{l-1}\sum_{m=1}^{k+1}
(-1)^{m}\bin{k}{m-1} m^{l-1} - \sum_{k=0}^{l-1} 2^{l-1-k}\sum_{m=1}^{k+1}
(-1)^{m}\bin{k}{m-1} m^{l-1}\\
& = &  2^{l-1}\sum_{m=1}^{l}
(-1)^{m}  m^{l-1}\sum_{k=m-1}^{l-1}\bin{k}{m-1} -
 \sum_{m=1}^{l}(-1)^{m}  m^{l-1} \sum_{k=m-1}^{l-1} 2^{l-1-k}
 \bin{k}{m-1}\\
& = &  2^{l-1}\sum_{m=1}^{l}
(-1)^{m}  m^{l-1}\bin{l}{m} -
 \sum_{m=1}^{l}(-1)^{m}  m^{l-1} \sum_{j=m}^{l}\bin{l}{j},\\
\end{eqnarray*}

since $\sum_{k=m-1}^{l-1}\bin{k}{m-1}=\bin{l}{m}$ and
$\sum_{k=m-1}^{l-1} 2^{l-1-k} \bin{k}{m-1}$ by Lemma~\ref{binrmq}.
By Lemma~\ref{vanlint}, 
$
\sum_{m=1}^{l}
(-1)^{m}  m^{l-1}\bin{l}{m} = 0
$ thus 
$$
\chi_{l,l}= - \sum_{m=1}^{l}(-1)^{m}  m^{l-1} \sum_{j=m}^{l}\bin{l}{j}\,
$$

and $\chi_{l,l} = (-1)^l \sigma_{l,l}$.
This is the desired
equality for $l$ even.
For an odd $l > 2$
we use the symmetry of the expression of
$\sigma_{l,l}$ and write
$$\displaylines{\sigma_{l,l} =
\sum_{p=0}^{\frac{l-3}{2}} \Big({(-1)}^{p+1} \sum_{q=0}^{p+1}
(-1)^q \bin{l}{q}(p+1-q)^{l -1} \cr
+ (-1)^{l-1-p} \sum_{r=0}^{l-1-p} (-1)^r
\bin{l}{r}(l-1-p-r)^{l-1}\Big).}$$
 With the change of index $q=l-r$ in the second term of the sum and
noticing that the contribution for $q=p+1$ is zero, one gets
 $$\sigma_{l,l} =
\sum_{p=0}^{\frac{l-3}{2}} {(-1)}^{p+1}
\sum_{q=0}^{l}(-1)^q\bin{l}{q} (p+1-q)^{l-1}.$$
The right hand side of the last equality is zero by
Lemma \ref{vanlint},
and thus, $\sigma_{l,l} = \chi_{l,l} = 0$. \CQFD

According to Lemma \ref{sigeqchi} the coefficients of $\sigma(Z)$
and $\chi(\RR Z)$ in the expressions of Lemma~\ref{exprsigchi} are
equal, and thus, $\sigma(Z) = \chi(\RR Z)$.\CQFD

\section{Nakajima polytopes and M-surfaces}\label{Mhypex}

In this section  we give examples of  families of $M$-surfaces
obtained by $T$-construction as hypersurfaces of $3$-dimensional toric
varieties.

\subsection{Newton polytopes without maximal hypersurfaces}\label{Mhypconterex}

We first explain that Itenberg-Viro's theorem (see \cite{IteVir}) on
the existence of $M$-hypersurfaces of any degree in the projective
spaces of any dimension cannot be generalized straightforwardly to all
projective toric varieties. Indeed, in dimensions greater or equal to $3$,
there exist polytopes $\Delta$ such that no hypersurface in $X_\Delta$
with Newton polytope $\Delta$ is maximal (See Proposition~\ref{nomax2}). Note that this does not mean
that the toric variety $X_\Delta$ does not admit $M$-hypersurfaces.

Clearly, if $\Delta$ is an interval $[a, b]$ in $\RR$, where $a$ and
$b$ are nonnegative integers, then there exists a maximal
$0$-dimensional subvariety in $\CC P^1 = X_\Delta$ with the Newton
polygon $\Delta$.

If $\Delta$ is a polygon in the first quadrant of $\RR^2$, then again
there exists a maximal curve in $X_\Delta$ with the Newton polygon
$\Delta$.  Such a curve can be constructed by the combinatorial
patchworking: it suffices to take as initial data a primitive convex
triangulation of $\Delta$ equipped with the following distribution of
signs: an integer point $(i, j)$ of $\Delta$ gets the sign ``-'' if
$i$ and $j$ are both even, and gets the sign ``+'', otherwise (cf.,
for example, \cite{Iterag}, \cite {IteVir2}, \cite{Haa}).

However, in dimension $3$ there are polytopes $\Delta$ such that no
surface in $X_\Delta$ with the Newton polytope $\Delta$ is
maximal (See Proposition~\ref{nomax}).

Let $k$ be a positive integer number, and $\Delta_k$ be the
tetrahedron in $\RR^3$ with vertices $(0,0,0), (1,0,0), (0,1,0)$, and
$(1,1,k)$.  Note that the only integer points of $\Delta_k$ are its
vertices.

\begin{figure}[htbp]
\begin{center}
\resizebox {3 cm}{!}{
\input{dessins/tetra.pstex_t}
}
\caption{Tetrahedron $\Delta_3$.}
\label{figure:fig1tetra}
\end{center}
\end{figure}
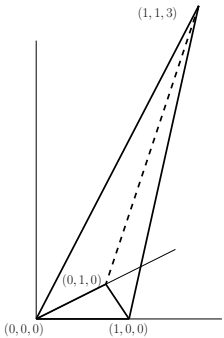

\begin{bbprop}\label{nomax}
  For any odd $k \geq 3$ and any even $k \geq 8$, there is no maximal
  surface in $X_{\Delta_k}$ with the Newton polytope $\Delta_k$.
\end{bbprop}

The proof of the above Proposition can be found in \cite{these} and
\cite{Moi1}. There, this family of examples is generalized to higher
dimensional polytopes and the following proposition is proved.

\begin{bbprop}\label{nomax2}
  For any integer $d \ge 3$ there exist $d$-dimensional polytopes
  $\Delta$ such that no hypersurface in $X_\Delta$ with the Newton
  polytope $\Delta$ is maximal.
\end{bbprop}

\subsection{Construction of $M$-surfaces}\label{S:msurf}

\begin{figure}[htbp]
\begin{center}
\resizebox {6 cm}{!}{
\input{dessins/nakapol2.pstex_t}
\hskip 8cm
\includegraphics{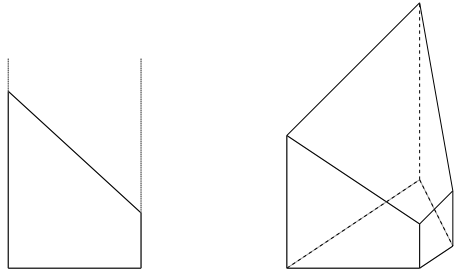}
}
\caption{Nakajima polytopes of dimension $2$ and $3$.}
\label{figure:nakapol3}
\end{center}
\end{figure}

Using Corollary~\ref{max} of Theorem~\ref{theprop} we prove the
existence of maximal surfaces in all nonsingular toric varieties
corresponding to Nakajima polytopes of dimension $3$ (see
Figure~\ref{figure:nakapol3}).


  A Nakajima polytope (See Introduction) is a {\bi
  nondegenerate} if it is $0$-dimensional or if $\bar{\Delta}$ is
nondegenerate and $f$ is positive on $\bar{\Delta}$.

\begin{figure}[htbp]
\begin{center}
\resizebox {6 cm}{!}{
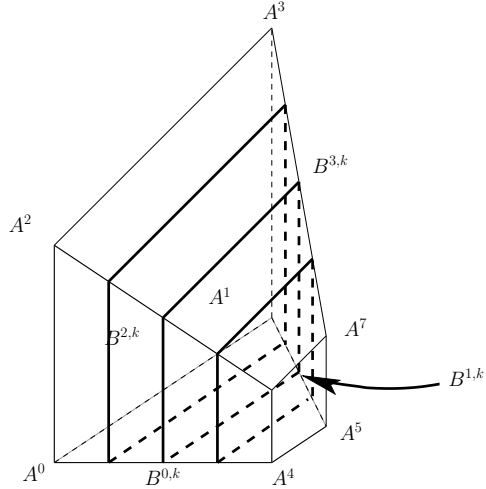
}
\caption{Slicing a nondegenerate Nakajima polytope.}
\label{figure:nakapol32}
\end{center}
\end{figure}

\begin{bbthm}\label{P2}
Let $\Delta$ be a $3$-dimensional Nakajima polytope corresponding to a
nonsingular toric variety $X_\Delta$.
Then, there exists a maximal surface
in $X_\Delta$ with the Newton polytope $\Delta$.
\end{bbthm}

\proof
Let  us first  assume that  $\Delta$  is
nondegenerate.

Decompose $\Delta$ into slices $T_k = \Delta \cap \{ (x,y,z)\in
(\RR^+)^3 ,k-1 \le x \le k \}$.  Let $s_k$ be the section $\Delta \cap
\{ (x,y,z)\in (\RR^+)^3, x = k\}$.  Denote by $B^{0,k}$, $B^{1,k}$,
$B^{2,k}$, and $B^{3,k}$ the vertices of $s_k$ (see Figure
\ref{figure:nakapol32}). We subdivide the slices $T_k$ in two cones
and two joins.  Take the cones of apex $B^{0,2k}$ and $B^{0,2k+2}$ on
$s_{2k+1}$ and the cones of apex $B^{3,2k-1}$ and $B^{3,2k+1}$ on
$s_{2k}$.  Take primitive convex triangulations of the sections $s_{k}$. They
induce a primitive triangulation of the cones.  The joins
$[B^{0,2k},B^{1,2k}] * [B^{1,2k+1}, B^{3,2k+1}]$, $[B^{0,2k},B^{2,2k}]
* [B^{2,2k+1}, B^{3,2k+1}]$, $[B^{0,2k+2},B^{1,2k+2}] * [B^{1,2k+1},
B^{3,2k+1}]$, $[B^{0,2k+2},B^{2,2k+2}] * [B^{2,2k+1}, B^{3,2k+1}]$ are
also naturally primitively triangulated which gives a primitive
convex triangulation of $\Delta$. Take the following
distribution of signs at the integer points of $\Delta$:

{\it any integer point $(i,j,k)$  gets ``-'' if $j$ and $k$
are both odd, and it gets ``+'', otherwise.}

The polytope $\Delta$ is now equipped with a
convex primitive triangulation and a sign distribution,
so  we  can apply  the  $T$-construction  as  in Section \ref{tcon}.

Let $p$ be an integer interior point of $\Delta$ and $\Star(p)$ be its
star. Assume that $p$ belongs to a section $s_{l}$. Then, $\Star(p)$
has two vertices $c_1$ and $c_2$ outside $s_{l}$. They are the apices
of the two cones over $s_{l}$.  They have the same parity (i.e. their
coordinates have the same reduction modulo $2$) and the distribution
of signs depends only on the parity thus, in each octant, their
symmetric copies carry the same sign.  Consider an octant where the
symmetric copy $q_0$ of $p$ is isolated in $s_{l}$ (i.e.  all vertices
of $\Star(q_0) \cap s_{l}$ except $q_0$ carry the sign opposite to the
sign of $q_0$).  Then, either in this octant the symmetric copies of
$c_1$ and $c_2$ carry the sign opposite to the sign of $q_0$ (and,
hence, $q_0$ is isolated), or $r(q_0)$ is isolated, where $r$ is the
reflection with respect to the coordinate plane $x = 0$.  Thus, for
each integer interior point $p$ of $\Delta$, there exists a symmetric
copy $q$ of $p$ such that $q$ is surrounded by a sphere $S^2(p) =
\Star(q) \cap H^*$ (see Remark~\ref{sphere}).  Moreover 
 at least one component of $H^*$ intersects
the coordinate planes.  Thus, the $T$-surface constructed has at least
$l^*(\Delta) +1$ connected components, and Corollary~\ref{max} shows
that this surface is maximal.
\begin{figure}[htbp]
\begin{center}
\resizebox {4 cm}{!}{
\input{dessins/star.pstex_t}
}
\caption{The star of $p$.}
\label{figure:star}
\end{center}
\end{figure}
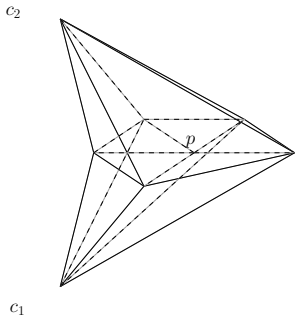

\begin{figure}[htbp]
\begin{center}
\resizebox {4 cm}{!}{
\input{dessins/nakatri2.pstex_t}
}
\caption{A degenerate Nakajima polytope.}
\label{figure:nakatri2}
\end{center}
\end{figure}
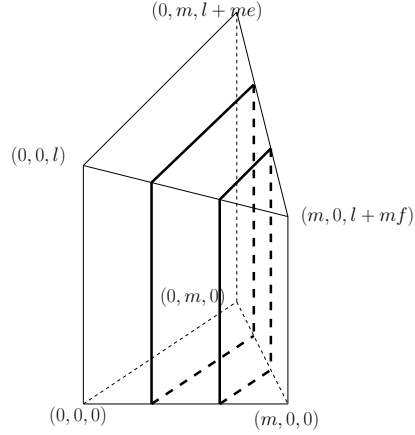

The degenerate case splits into two subcases (Since we only
consider here Nakajima polytopes corresponding to nonsingular toric
varieties).
  Either the Nakajima
polytope is a truncated cylinder over a triangle corresponding to the
projective plane, or it is a tetrahedron corresponding to the
projective $3$-space.  The existence of $M$-surfaces in the latter
case was proved by Viro in \cite{Vir79}. In the former case the
Nakajima polytope $\Delta$ is the convex hull of the triangles
$((0,0,0), (m,0,0), (0,m,0))$ and $((0,0,l), (0,m,l+me), (m,0,l+mf))$
for some integers $m$, $l$, $e$ and $f$ (see Figure
\ref{figure:nakatri2}).

Decompose $\Delta$  into slices $T_k =
\Delta \cap \{ (x,y,z)\in (\RR^+)^3 ,k-1 \le x \le k \}$. Let $s_k$ be
the section $\Delta \cap \{ (x,y,z)\in (\RR^+)^3, x = k\}$.
We triangulate  $\Delta \cap \{ (x,y,z)\in  (\RR^+)^3 ,y \le  m-1 \}$ using
the same triangulation as in the nondegenerate case. Take the
cone over $s_{m-1}$  of apex $(m,0,0)$ if $m$  is even or $(m,0,l+mf)$,
otherwise. Take in the cone the triangulation induced by the
triangulation of $s_{m-1}$.
Subdivide the only remaining non-primitive tetrahedron
into primitive ones in
the unique possible way.
Take the same distribution of
signs that in the nondegenerate case.

Then, as  in the nondegenerate case, any  interior point has an
isolated symmetric  copy. There is  also a component  intersecting the
coordinate planes. Thus, the surface is maximal.
\CQFD

Note that Theorem \ref{P2}
produces,
in particular,
$M$-surfaces
in $\Sigma_\alpha \times \CC P^1$.

For a nonnegative integer $\alpha$ and positive integers $l$, $m$ and
$n$ let $\delta_\alpha^{m,n}$ and $\Delta^{l,m,n}_\alpha$ be the
polytopes defined in Section~\ref{S:tropstate}. Recall that the toric
variety $X_{\Delta^{l,m,n}_\alpha}$ is isomorphic to $\Sigma_\alpha
\times \CC P^1$.  Since $\Delta^{l,m,n}_\alpha$ is a Nakajima
polytope, the following statement is a corollary of Theorem \ref{P2}.

\begin{bbcor}\label{hirz}

For  any
nonnegative integer $\alpha$ and any positive integers
$m$, $n$ and $l$,
there exists a maximal
surface 
in $\Sigma_\alpha \times \CC P^1$
with the Newton polytope $\Delta^{l,m,n}_\alpha$.\CQFD
\end{bbcor}

\section{Appendix}\label{appendix}

\subsection{Usefull combinatorial lemmae}

We insert here the proofs of lemmae that we use quite often, mainly
in Section~\ref{sigmachi}.

\begin{bblem}\label{vanlint}
Let $l$ and  $i$ be
nonnegative integers. Then,
for $l + 1 \le i$, one has (See \cite{LiW} p. 71)
\begin{equation}\label{vl}
\sum_{q=0}^{i} (-1)^q \bin{i}{q} q^{l} \;= \; \sum_{q=0}^{i} (-1)^q
\bin{i}{q}(i-q)^{l} = 0
\end{equation}
and, as a consequence, for any integer $p$,
\begin{equation}\label{nbsurj}
\sum_{q=0}^{i} (-1)^q
\bin{i}{q}(p-q)^{l} = 0.
\end{equation}
\end{bblem}

\proof
Formula \ref{vl} can be found in  \cite{LiW} p. 71.
We prove Equality (\ref{nbsurj}) using (\ref{vl}).

Write $p-q  =  (i-q)+(p-i)$ and then,
\begin{eqnarray*}
(p-q)^l & = & \sum_{m=0}^{l} \bin{l}{m}
(i-q)^{l-m}(p-i)^m .\\
\end{eqnarray*}
So,
\begin{eqnarray*}
\sum_{q=0}^{i} (-1)^q \bin{i}{q}(p-q)^{l}& = &  \sum_{q=0}^{i} (-1)^q
\bin{i}{q} \sum_{m=0}^{l} \bin{l}{m}
(i-q)^{l-m}(p-i)^m\\
& =  & \sum_{m=0}^{l} \bin{l}{m}(p-i)^m \Bigg(\sum_{q=0}^{i}
(-1)^q \bin{i}{q} (i-q)^{l-m}
\Bigg)
\end{eqnarray*}
The last term is zero by Equality \ref{vl}. \CQFD

{\bf Proof of lemma \ref{psi}}. 
From the definitions we see that $$Q_\Delta(t) = (\sum_{k=0}^{d+1}
\bin{d+1}{k} t^k (-1)^k) \cdot \sum_{n=0}^\infty (\sum_{i=0}^d a_i n^i)
t^n,$$
and thus, $\Psi_j = \sum_{i=0}^d (\sum_{n=0}^j
(-1)^{j-n}\bin{d+1}{j-n} n^i) a_i^\Delta$.
Let $A_j = \sum_{n=0}^j
(-1)^{j-n}\bin{d+1}{j-n} n^i$.\\
If $j \ge d+1$, then
$$ A_j=\sum_{k=0}^{d+1}
(-1)^k \bin{d+1}{k} (j-k)^i,$$
which is zero according to
Lemma~\ref{vanlint} Equality~\ref{nbsurj}. \CQFD

\begin{bblem}\label{binrmq}
One has $\sum_{n=0}^{p} 2^{p-n}
\bin{n}{k} = \sum_{l=k+1}^{p+1} \bin{p+1}{l}$.
\end{bblem}

\proof
We use the fact that $\sum_{n=k}^{p} \bin{n}{k} = \bin{p+1}{k+1}$
and write
\begin{eqnarray*}
\sum_{n=0}^{p} 2^{p-n}\bin{n}{k} & = & \sum_{n=k}^{p} 2^{p-n}\bin{n}{k}\\
& = & \sum_{n=k}^{p} \bin{n}{k} + \sum_{i=0}^{p} (2^i)\sum_{n=k}^{p-i}\bin{n}{k}\\
& = & \bin{p+1}{k+1} + \sum_{i=0}^{p-k}(2^i) \sum_{i=0}^{p-k} \bin{p+1-i}{k+1}\\
& = & \bin{p+1}{k+1} + \sum_{n=k+1}^{p+1} 2^{p+1-n}\bin{n}{k+1},\\
\end{eqnarray*}
and the result follows by induction.\CQFD

\subsection{Formula for the number of $k$-simplices}\label{ksimplices}

In order to prove the formula of Proposition~\ref{nbsrformula} we
first check that it is true for all primitive simplices.

\subsubsection {Formula for a primitive simplex}

Let $s_i$ be a primitive simplex of
dimension $i$.
Then
$$Ehr_{s_i}(\lambda) = \frac{1}{i !}(\lambda + 1) \dots
(\lambda + i)=\frac{1}{i !} \sum_{j=1}^{i+1}{(-1)}^{i+1-j}S_1(i+1,j)
\lambda^{j-1}, $$
where $S_1$ is the first Stirling number defined by the formula
$\sum_{m=0}^n S_1(n,m)x^m= x (x-1) \dots (x-n+1)$ (see \cite{BV}).
Thus
$$a^{s_i}_{j-1}=\frac{(-1)^{i+1-j}}{i !}S_1(i+1,j).$$
It remains to show that
$$ N_{r}^{s_i}=\frac{r !}{i!}
\sum_{l=r+1}^{i+1} S_2(l,r+1).S_1(i+1,l).$$
We have
$\sum_{l=r+1}^{i+1} S_2(l,r+1).S_1(i+1,l)=\delta_{r,i}$ where
$\delta_{r,i}$ is the Kronecker index (see
\cite{LiW} p. 107), and thus,
\[\frac{r !}{i!}\sum_{l=r+1}^{i+1} S_2(l,r+1).S_1(i+1,l)
=\delta_{r,i}\]
 which is exactly the number of $r$-dimensional
simplices contained in the interior of $s_i$ (zero if $r \neq i$ and one othewise).

\subsubsection{Proof for a general polytope}

\begin{bbdef}
  A triangulation $\tau$ of a polytope of dimension $d$ is called {\bi
    shellable} if there exists a numbering of its $d$-simplices $s_1,
  s_2 , \dots ,s_k$ such that for $i \in \{2, \dots, k\}$
 $$ s_i \cap \cup_{j=1}^{i-1} s_j$$
is a
nonempty union of $(d-1)$-simplices of $\tau$ homeomorphic to a
$(d-1)$-ball. This numbering is called a {\it shelling} of $\tau$.
\end{bbdef}

\begin{bbthm}[Ziegler (\cite{zie}, p. 243)]
Every convex triangulation of a
polytope is shellable.
\end{bbthm}

The formula of Proposition~\ref{nbsrformula} holds for a point (with
the convention that the point is in its interior).  We now assume that
the formula is true in all dimensions less than $d$.  Let $\Delta$ be
a $d$-dimensional polytope, and $\tau$ be a primitive triangulation of
$\Delta$.  Fix $s_1,s_2,\dots, s_t$ a shelling of $\tau$, and put
$U_j=\cup_{i=1}^j s_j$.  Note that $U_j$ need not to be convex.
Assume that the formula is true for $U_j$.  Note that the formula also
holds for $U_j\cap s_{j+1}$ which is $(d-1)$-dimensional, eventhough
it is not a polytope.  We know that the formula is also true for
$s_{j+1}$.

\begin{bblem}
The numbers $N_r$
satisfy the relation
$$N_r^{U_j} + N_r^{s_{j+1}} + N_r^{U_j
\cap s_{j+1}} = N_r^{U_{j+1}}.$$
\end{bblem}

\proof

This follows immediately from the fact that
$   \stackrel{\circ}{U_j}   \sqcup  \stackrel{\circ}{s_{j+1}}   \sqcup
\stackrel{\circ}{  ( U_j  \cap s_{j+1})}  = \stackrel{\circ}{U_{j+1}}$,
where $\stackrel{\circ}{U}$ stands for the interior of $U$.
\CQFD

By induction hypothesis,
\begin{eqnarray*}
 N_r^{U_j \cap s_{j+1}} & = & \sum_{l=r+1}^{d} \stir(r+1,l)
 (-1)^{d-l} a_{l-1}^{U_j \cap s_{j+1}},\\
\end{eqnarray*}
and
since $ a_{d}^{U_j \cap s_{j+1}}=0$, we can also write
\begin{eqnarray*}
 N_r^{U_j \cap s_{j+1}} & = & \sum_{l=r+1}^{d+1} \stir(r+1,l)
 (-1)^{d-l} a_{l-1}^{U_j \cap s_{j+1}}.\\
\end{eqnarray*}

Then $N_r^{U_{j+1}} =  \sum_{l=r+1}^{d+1}
\stir(r+1,l)(-1)^{d-l+1}.(a^{U_j}_{l-1}+ a^{s_{j+1}}_{l-1}- a^{U_j \cap
s_{j+1}}_{l-1})$,
and $a^{U_j}_{l-1}+ a^{s_{j+1}}_{l-1}- a^{U_j \cap
s_{j+1}}_{l-1}$ is precisely the coefficient $a^{U_{j+1}}_{l-1}$ in the
Ehrhart polynomial of $U_{j+1}$ (see \cite{Ehr67} for Ehrhart
polynomials of general polyhedra). This
completes the proof of the formula. \proofend

\def\cprime{$'$} \def\cprime{$'$} \def\cprime{$'$}


\end{document}

%% file: dessins/delta.pstex_t
\begin{picture}(0,0)%
\includegraphics{dessins/delta.pstex}%
\end{picture}%
\setlength{\unitlength}{3947sp}%
\begingroup\makeatletter\ifx\SetFigFont\undefined%
\gdef\SetFigFont#1#2#3#4#5{%
  \reset@font\fontsize{#1}{#2pt}%
  \fontfamily{#3}\fontseries{#4}\fontshape{#5}%
  \selectfont}%
\fi\endgroup%
\begin{picture}(7662,4596)(751,-5719)
\put(751,-5536){\makebox(0,0)[lb]{\smash{\SetFigFont{20}{24.0}{\familydefault}{\mddefault}{\updefault}$(0,0)$}}}
\put(751,-1411){\makebox(0,0)[lb]{\smash{\SetFigFont{20}{24.0}{\familydefault}{\mddefault}{\updefault}$(0,m)$}}}
\put(3226,-1411){\makebox(0,0)[lb]{\smash{\SetFigFont{20}{24.0}{\familydefault}{\mddefault}{\updefault}$(n,m)$}}}
\put(8401,-5611){\makebox(0,0)[lb]{\smash{\SetFigFont{20}{24.0}{\familydefault}{\mddefault}{\updefault}$(n+m\alpha,0)$}}}
\put(2326,-3586){\makebox(0,0)[lb]{\smash{\SetFigFont{29}{34.8}{\familydefault}{\mddefault}{\updefault}$\delta_\alpha^{m,n}$}}}
\end{picture}

%% file: dessins/P.pstex_t
\begin{picture}(0,0)%
\includegraphics{dessins/P.pstex}%
\end{picture}%
\setlength{\unitlength}{3947sp}%
\begingroup\makeatletter\ifx\SetFigFont\undefined%
\gdef\SetFigFont#1#2#3#4#5{%
  \reset@font\fontsize{#1}{#2pt}%
  \fontfamily{#3}\fontseries{#4}\fontshape{#5}%
  \selectfont}%
\fi\endgroup%
\begin{picture}(9537,8295)(376,-8044)
\put(9901,-7411){\makebox(0,0)[lb]{\smash{\SetFigFont{20}{24.0}{\familydefault}{\mddefault}{\updefault}$(n+m\alpha,0,0)$}}}
\put(451,-7936){\makebox(0,0)[lb]{\smash{\SetFigFont{20}{24.0}{\familydefault}{\mddefault}{\updefault}$(0,0,0)$}}}
\put(2401,-5761){\makebox(0,0)[lb]{\smash{\SetFigFont{20}{24.0}{\familydefault}{\mddefault}{\updefault}$(0,m,0)$}}}
\put(6601,-5761){\makebox(0,0)[lb]{\smash{\SetFigFont{20}{24.0}{\familydefault}{\mddefault}{\updefault}$(n,m,0)$}}}
\put(376,-2311){\makebox(0,0)[lb]{\smash{\SetFigFont{20}{24.0}{\familydefault}{\mddefault}{\updefault}$(0,0,l)$}}}
\put(2476,-736){\makebox(0,0)[lb]{\smash{\SetFigFont{20}{24.0}{\familydefault}{\mddefault}{\updefault}$(0,m,l)$}}}
\put(6526,-736){\makebox(0,0)[lb]{\smash{\SetFigFont{20}{24.0}{\familydefault}{\mddefault}{\updefault}$(n,m,l)$}}}
\put(9751,-2686){\makebox(0,0)[lb]{\smash{\SetFigFont{20}{24.0}{\familydefault}{\mddefault}{\updefault}$(n+m\alpha,0,l)$}}}
\put(3976,-4036){\makebox(0,0)[lb]{\smash{\SetFigFont{29}{34.8}{\familydefault}{\mddefault}{\updefault}$P_\alpha^{l,m,n}$}}}
\end{picture}

%% file: dessins/tetra.pstex_t
\begin{picture}(0,0)%
\includegraphics{dessins/tetra.pstex}%
\end{picture}%
\setlength{\unitlength}{3947sp}%
\begingroup\makeatletter\ifx\SetFigFont\undefined%
\gdef\SetFigFont#1#2#3#4#5{%
  \reset@font\fontsize{#1}{#2pt}%
  \fontfamily{#3}\fontseries{#4}\fontshape{#5}%
  \selectfont}%
\fi\endgroup%
\begin{picture}(5637,8627)(376,-8044)
\put(376,-7936){\makebox(0,0)[lb]{\smash{\SetFigFont{20}{24.0}{\familydefault}{\mddefault}{\updefault}$(0,0,0)$}}}
\put(3076,-7936){\makebox(0,0)[lb]{\smash{\SetFigFont{20}{24.0}{\familydefault}{\mddefault}{\updefault}$(1,0,0)$}}}
\put(1876,-6661){\makebox(0,0)[lb]{\smash{\SetFigFont{20}{24.0}{\familydefault}{\mddefault}{\updefault}$(0,1,0)$}}}
\put(3826,239){\makebox(0,0)[lb]{\smash{\SetFigFont{20}{24.0}{\familydefault}{\mddefault}{\updefault}$(1,1,3)$}}}
\end{picture}

%% file: dessins/nakapol2.pstex_t
\begin{picture}(0,0)%
\includegraphics{dessins/nakapol2.pstex}%
\end{picture}%
\setlength{\unitlength}{3947sp}%
\begingroup\makeatletter\ifx\SetFigFont\undefined%
\gdef\SetFigFont#1#2#3#4#5{%
  \reset@font\fontsize{#1}{#2pt}%
  \fontfamily{#3}\fontseries{#4}\fontshape{#5}%
  \selectfont}%
\fi\endgroup%
\begin{picture}(3644,5734)(2379,-5183)
\end{picture}

%% file: dessins/decompnaka.pstex_t
\begin{picture}(0,0)%
\includegraphics{dessins/decompnaka.pstex}%
\end{picture}%
\setlength{\unitlength}{3947sp}%
\begingroup\makeatletter\ifx\SetFigFont\undefined%
\gdef\SetFigFont#1#2#3#4#5{%
  \reset@font\fontsize{#1}{#2pt}%
  \fontfamily{#3}\fontseries{#4}\fontshape{#5}%
  \selectfont}%
\fi\endgroup%
\begin{picture}(7275,8085)(1351,-8044)
\put(1576,-7861){\makebox(0,0)[lb]{\smash{\SetFigFont{20}{24.0}{\familydefault}{\mddefault}{\updefault}$A^0$}}}
\put(3601,-7936){\makebox(0,0)[lb]{\smash{\SetFigFont{20}{24.0}{\familydefault}{\mddefault}{\updefault}$B^{0,k}$}}}
\put(5701,-7936){\makebox(0,0)[lb]{\smash{\SetFigFont{20}{24.0}{\familydefault}{\mddefault}{\updefault}$A^4$}}}
\put(6826,-7186){\makebox(0,0)[lb]{\smash{\SetFigFont{20}{24.0}{\familydefault}{\mddefault}{\updefault}$A^5$}}}
\put(8626,-6286){\makebox(0,0)[lb]{\smash{\SetFigFont{20}{24.0}{\familydefault}{\mddefault}{\updefault}$B^{1,k}$}}}
\put(6901,-5461){\makebox(0,0)[lb]{\smash{\SetFigFont{20}{24.0}{\familydefault}{\mddefault}{\updefault}$A^7$}}}
\put(6376,-2761){\makebox(0,0)[lb]{\smash{\SetFigFont{20}{24.0}{\familydefault}{\mddefault}{\updefault}$B^{3,k}$}}}
\put(5551,-211){\makebox(0,0)[lb]{\smash{\SetFigFont{20}{24.0}{\familydefault}{\mddefault}{\updefault}$A^3$}}}
\put(4651,-4936){\makebox(0,0)[lb]{\smash{\SetFigFont{20}{24.0}{\familydefault}{\mddefault}{\updefault}$A^1$}}}
\put(2926,-5611){\makebox(0,0)[lb]{\smash{\SetFigFont{20}{24.0}{\familydefault}{\mddefault}{\updefault}$B^{2,k}$}}}
\put(1351,-3736){\makebox(0,0)[lb]{\smash{\SetFigFont{20}{24.0}{\familydefault}{\mddefault}{\updefault}$A^2$}}}
\end{picture}

%% file: dessins/star.pstex_t
\begin{picture}(0,0)%
\includegraphics{dessins/star.pstex}%
\end{picture}%
\setlength{\unitlength}{3947sp}%
\begingroup\makeatletter\ifx\SetFigFont\undefined%
\gdef\SetFigFont#1#2#3#4#5{%
  \reset@font\fontsize{#1}{#2pt}%
  \fontfamily{#3}\fontseries{#4}\fontshape{#5}%
  \selectfont}%
\fi\endgroup%
\begin{picture}(5197,5721)(1726,-6919)
\put(1726,-1486){\makebox(0,0)[lb]{\smash{\SetFigFont{20}{24.0}{\familydefault}{\mddefault}{\updefault}$c_2$}}}
\put(1801,-6811){\makebox(0,0)[lb]{\smash{\SetFigFont{20}{24.0}{\familydefault}{\mddefault}{\updefault}$c_1$}}}
\put(4951,-3811){\makebox(0,0)[lb]{\smash{\SetFigFont{20}{24.0}{\familydefault}{\mddefault}{\updefault}$p$}}}
\end{picture}

%% file: dessins/nakatri2.pstex_t
\begin{picture}(0,0)%
\includegraphics{dessins/nakatri2.pstex}%
\end{picture}%
\setlength{\unitlength}{3947sp}%
\begingroup\makeatletter\ifx\SetFigFont\undefined%
\gdef\SetFigFont#1#2#3#4#5{%
  \reset@font\fontsize{#1}{#2pt}%
  \fontfamily{#3}\fontseries{#4}\fontshape{#5}%
  \selectfont}%
\fi\endgroup%
\begin{picture}(5100,7596)(1126,-6844)
\put(1801,-6661){\makebox(0,0)[lb]{\smash{\SetFigFont{20}{24.0}{\familydefault}{\mddefault}{\updefault}$(0,0,0)$}}}
\put(5401,-6736){\makebox(0,0)[lb]{\smash{\SetFigFont{20}{24.0}{\familydefault}{\mddefault}{\updefault}$(m,0,0)$}}}
\put(3751,-4561){\makebox(0,0)[lb]{\smash{\SetFigFont{20}{24.0}{\familydefault}{\mddefault}{\updefault}$(0,m,0)$}}}
\put(3601,464){\makebox(0,0)[lb]{\smash{\SetFigFont{20}{24.0}{\familydefault}{\mddefault}{\updefault}$(0,m,l+me)$}}}
\put(1126,-2086){\makebox(0,0)[lb]{\smash{\SetFigFont{20}{24.0}{\familydefault}{\mddefault}{\updefault}$(0,0,l)$}}}
\put(6226,-3136){\makebox(0,0)[lb]{\smash{\SetFigFont{20}{24.0}{\familydefault}{\mddefault}{\updefault}$(m,0,l+mf)$}}}
\end{picture}